\documentclass[10pt]{article}

\usepackage{latexsym,color,amsmath,amsthm,amssymb,amscd,amsfonts}

\newtheorem{theorem}{Theorem}

\newtheorem{lemma}{Lemma}

\newtheorem{proposition}{Proposition}
\newtheorem{definition}{Definition}
\newtheorem{remark}{Remark}
\newtheorem{example}{Example}

\def\P{\mathfrak{P}_n}
\def\A{\mathfrak{A}_n}
\def\R{{\mathbb R}}
\def\K{{{\mathcal K}_n}}

\def\inte{{\rm int}}
\def\eps{{\varepsilon}}

\def\conv{{\rm conv}}
\def\vol{{\rm vol}}

\begin{document}

\title{Dual Affine invariant points
}

\author{Mathieu Meyer, Carsten Sch\"utt  and Elisabeth M. Werner 
\thanks{Partially supported by an NSF grant}}

\date{}

\maketitle

\date{}

 \maketitle

\begin{abstract}
An affine invariant point  on the class of convex bodies $\K$ in $R^n$, endowed with the Hausdorff metric,  is a continuous map
from $\K$ to $\R^n$  which is invariant under  one-to-one affine transformations $A$ on $ \R^n$, that is,  $p\big(A(K)\big)=A\big(p(K)\big)$.
\par
We define here the  new notion of dual affine point $q$ of an affine invariant point $p$ by  the  formula $q(K^{p(K)})=p(K)$
for every $K\in \K$, where $K^{p(K)}$ denotes the polar of $K$ with respect to $p(K)$.
\par
We investigate  which affine invariant points do have a dual point, whether this dual point is unique and has itself a dual point.
We define a product on the set of affine invariant points, in relation with duality.  
\par
Finally, examples are given which exhibit the rich structure of the set of affine invariant points.
\end{abstract}

\newpage
\section{ Introduction.}
\par
While convex bodies have 
been  the topic  of  extensive research for more than
a century, it is the affine geometry of these bodies that has been a main  focus of study in recent years.
We only mention the rapid progress in the $L_p$ Brunn Minkowski theory (e.g.,  \cite{Ga3, Hab, Lu2, LYZ2004, SW5, SA2, WernerYe2008, WernerYe2010}) and the theory of valuations e.g., \cite{Lud2, Lud3, LR2,Schuster2010}.
The resulting body of work has proved to be a
valuable tool in fields such as harmonic analysis, information theory, stochastic geometry and PDEs (e.g., \cite{LYZ2002/1, NPRZ2010, PaourisWerner2010, Werner2012}).
\par
Specific points associated to a convex body, like the centroid and
the Santal\'o point, satisfy an affine invariance property: The point
of an affine image
of a convex body is the affine image of the point. More formally,
 if $\K$ denotes the set of all convex bodies in $\R^n$, a mapping $p:\K\to \R^n$ 
 is {\it an affine invariant point} if $p$ is continuous  for the Hausdorff topology and satisfies 
$$p(T(K))=T\big(p(K)\big)$$
 for every one-to-one affine mapping $T:\R^n\to \R^n$.
 \par
Even though  this notion is intriguing in its simplicity, little  is known about affine invariant points.  
At the same time,  these are fundamental invariants
of convex sets.  They are, for instance, useful to characterize properties of symmetry or 
of non symmetry of convex bodies (e.g., \cite{MeyerSchuettWerner2011}  and \cite{MeyerSchuettWerner2012}).
The more different affine invariant points a convex body has the less symmetric it is. The task of computing an affine invariant point of a convex body can be formidable,
even to show that two affine invariant points of a convex body are different can be nontrivial.
\par
Affine invariant points were first defined by 
B. Gr\"unbaum  in 1963 in his seminal paper \cite{Gruenbaum1963},  where he also posed several open problems.
 In two preceding papers, \cite{MeyerSchuettWerner2011}  and \cite{MeyerSchuettWerner2012},  we answered some of Gr\"unbaum's questions: The dimension of the space of
 affine invariant points is infinte and there are convex bodies $K$ in $\mathbb R^{n}$
 such that every point in $\mathbb R^{n}$ is an affine invariant point of $K$.
 More importantly, we showed in some cases that the presence of many affine invariant points
 means that the convex body lacks symmetry.
 \par
 However,  many structural questions are still open.
 In this paper we address  them through the study of duality. We introduce the new notion of  {\em dual affine invariant point}. In short, 
the point  $q$ is a dual affine invariant point to $p$ if 
 $$
 q(K^{p(K)})=p(K)
 $$
 for all convex bodies $K$. Here,  
 $
 K^{z}= \{y\in \R^n: \langle y-z, x-z)\rangle \le 1\hbox{ for every }x\in K\}
 $
 is the polar of $K$ with respect to the point $z$.
The motivation for our definition, given in Section \ref{Section:dual},  comes from the duality relation between the center of gravity and the Santal\'o point of a convex body.  Further examples of dual affine invariant points are the center of the John ellipsoid and the center
of the L\"owner ellipsoid. All this is explained in Section \ref{example}.
\par
We start  our study by addressing a number of basic questions. First, 
does every affine invariant point  $p$ have  a dual $p^\circ$? The answer, surprisingly,  is: No. This is the content of Theorem \ref{theo:nichtinjective}. In Theorem \ref{theo:surjective},  we show that  if a dual affine invariant point   exists, it is unique.
Theorem \ref{theo:dual} establishes a reflexivity principle for affine invariant points, namely that the double dual $p^{\circ \circ}$ of $p$ equals $p$.  The proofs of the  theorems require a number of technical results. Those are  presented in Section \ref{Section:dual}.
\par
\par
In Section \ref{main} we give the proof of the main theorems. 
We also define  there a product $[p,q]$ of two affine invariant points $p$ and $q$ as a mapping from from the set $\mathcal{P}_n$ of all affine invariant points into itself.  This product has a nice duality property, 
$$
[q^{\circ},p^{\circ}]\circ [p,q]= I_n,
$$
where $I_n$ is the identity on $\mathcal{P}_n$.
\par
Finally, Section \ref{example} is devoted  to a list of useful examples of  affine invariant points and a related notion, that of  {\em affine invariant sets} (also defined in Section \ref{background}),   many new ones among them.
Moreover, we investigate how to extend 
the affine invariant points on 
${\cal K}_{n,k}$, the set of compact convex subsets in $\R^n$ whose affine span is $k$-dimensional, to affine invariant points on $\mathcal K_n$.
\vskip 4mm
The authors would like to thank the American Institute of Mathematics in Palo Alto where, in the course of  the workshop 
``Invariants in  convex geometry
and  Banach space theory", 
much of the paper was produced.

\section
{Notation and Background material.} \label{background}
\par
We  denote by $\mathcal K_n$  the set of all convex bodies in 
$\mathbb R^{n}$, that is the set of all  convex compact sets
with nonempty interior.  
For $K \in \mathcal K_n$,  $\mbox{int}(K)$ is the interior of $K$ and $\partial K$  is its boundary.
We say that $K  \in \mathcal{K}_n$ is in $C^2_+$,  if $\partial K$ is $C^2$  with strictly positive 
Gaussian curvature.
\par
The Euclidean ball centered at $a$ with radius $r$ is  $B^n_2(a,r)$.  We write in short  $B^n_2=B^n_2(0,1)$ and $S^{n-1}=\partial B^n_2$. We endow $\R^n$ with its canonical scalar product, and for $x\in \R^n$, we denote $|x|=\sqrt{\langle x,x\rangle}$ its Euclidean norm. 
The $n$-dimensional volume of $K$ is $\text{vol}_n(K)$, or simply $|K|$. Quite often, if $A\subset \R^n$ has an affine span
of dimension $k$,  we shall denote also by $|A|$ the $k$-dimensional volume of $A$ in its affine span.
\par
For subsets  $A$ and $B$ of $\R^n$, 
$\conv[A,B]$ denotes their convex hull, the smallest convex body containing them.
\par
The  support function $h_{K}:\R^{n}
\to \R$ of a  convex body $K$  is given by
$
h_{K}(\xi)=\sup_{x\in K} \ \langle \xi,x \rangle$. 
If $0\in \inte(K)$,
$K^\circ=\{y \in \mathbb R^{n}: \langle x,y \rangle \leq 1\}$ 
 is the polar body of $K$ with respect to $0$.
More generally,  we define the polar body $K^x$ of $K$ with respect to $x \in \mathbb{R}^n$ by 
$$K^x=(K-x)^\circ+x, \hbox{ \hskip 1mm or \hskip 1mm } K^x-x=(K-x)^\circ.$$
Note that  $K^x\in \K$ if and only if  $x \in \mbox{int}(K)$.  We  will only consider such situations. By the bipolar theorem, 
\begin{equation}\label{bipolar}
(K^x)^x=K,
\end{equation}
 which may be written as 
$$K-x=(K^x-x)^{\circ}.$$
\par
We shall frequently use  the  fact that  if $T:\R^n\to \R^n$ is a one-to-one linear map, $K\in \K$ and $x\in\inte(K)$, then 
\begin{equation}\label{adjoint}
\big(T(K-x)\big)^{\circ}= T^{*-1}\big((K-x)^{\circ}\big).
\end{equation}
Here $T^*$ is the adjoint  of $T$ and $T^{*-1}$  its inverse. 
\par
The Hausdorff metric $d_H$ on $\mathcal K_n$ is defined as
\begin{equation*}\label{Hausdorff}
d_H(K_1, K_2) = \min\{\lambda \geq 0: K_1 \subseteq K_2+\lambda B^n_2,    K_2 \subseteq K_1+\lambda B^n_2 \}.
\end{equation*}
\vskip 3mm
Now we  recall the definitions of affine invariant points and of affine invariant sets \cite{Gruenbaum1963, MeyerSchuettWerner2012}.
\vskip 2mm
\begin{definition}\label{aip}
A  map 
$p:{\mathcal K}_n \rightarrow\mathbb R^{n}$ is called affine invariant point,  if $p$ is continuous and if
for every nonsingular affine map $T:\mathbb R^{n}\rightarrow \mathbb R^{n}$, one has
\begin{equation*} \label{def1}
p\big(T(K)\big)=T\big(p(K)\big).
\end{equation*}
We denote by $\mathfrak{P}_n$ the set of affine invariant points in $\mathbb{R}^n$,
\begin{equation*} \label{def:aip}
\mathfrak{P}_n =\{ p: \mathcal K_n \rightarrow\mathbb R^{n}  \big| \  p \  \text{ is  continuous and affine invariant} \}, 
\end{equation*}
and for a fixed $K \in \K$, $\mathfrak{P}_n(K)=\{p(K);   p \in \mathfrak{P}_n\}$.
\par
We say that  $p\in \P$  is proper  if for all
$K\in {\mathcal K}_n$, one has $p(K)\in \operatorname{int}(K)$.
\end{definition}
\vskip 3mm
\begin{definition}\label{aism} A map $A:\K\to \K$ is an affine invariant set mapping, or an affine invariant set, if $A$ is continuous (when $\K$ is endowed with the Hausdorff metric) and if for every affine one-to-one map $T:\R^n\to \R^n$, one has
$$A\big(T(K)\big)= T\big(A(K)\big).$$
We denote by $\mathfrak{A}_n$ the set of all affine invariant set mappings from $\K$ to $\K$.
\end{definition}
\vskip 3mm
Well known classical examples (see e.g. \cite{Gruenbaum1963}, \cite{MeyerSchuettWerner2012}) of  proper   affine invariant points  of a convex body
$K$ in
$\mathbb R^{n}$ are the  {\em centroid}, 
\begin{equation}\label{centroid}
g(K)
=\frac{\int_{K}xdx}{|K|},
\end{equation}
the  {\em Santal\'o point}, which is 
the unique point $s(K)\in \inte(K)$ for which $|K^{s(K)}|=\min_{x} |K^x|$,  
  the center of the  {\em John ellipsoid} of $K$, that is the {\em ellipsoid of maximal volume}
contained in $K$ and 
the center of the  {\em L\"owner ellipsoid} of $K$, that is the {\em ellipsoid of minimal volume} containing $K$. 
\par
We will discuss these and other examples  in Section \ref{example}. 
More details on affine invariant points, and some results that we shall use here,   
can be found in \cite{MeyerSchuettWerner2012}.

\vskip 3mm

\section{Dual affine invariant points.} \label{Section:dual}

We now introduce the new concept of a dual of an affine invariant point.
\vskip 2mm

\begin{definition}\label{daip}
Let $p \in \mathfrak{P}_n$ be proper.
\vskip 1mm
\noindent
(i) We say $q\in \P$ is a dual of $p$
if for all  $K\in \K$
\begin{equation}\label{dual point1}
q(K^{p(K)})=p(K).
\end{equation}
If $p$ has a unique dual, we denote it by $p^{\circ}$.
\vskip 1mm
\noindent
(ii) Let $w\in \P$. We say that $w$ is a bidual of $p$,  if $p$ has a proper dual $q\in \P$, such that $w$ is a dual of $q$. This means that
there is an affine invariant  point $q$ that is dual to $p$ and that $w$ is dual to $q$, 
$$
q(K^{p(K)})=p(K) \ \  \   \text{and} \ \ \  w(K^{q(K)})=q(K).
$$
\end{definition}
\vskip 2mm
Theorem \ref{theo:dual} assures that a dual point is automatically proper.
\par
The centroid and  the  Santal\'o point,  and  the center of the  John ellipsoid of $K$ and 
the center of the   L\"owner ellipsoid  are examples of dual affine invariant points.  We will explain this  in Section \ref{example}.
\vskip 3mm
For   $q\in \P$, the affine invariance of the 
mapping $K\to q(K^{p(K)})$ implies that $q$ is a dual of $p$. This is the content of the following lemma.
\vskip 2mm
\begin{lemma} Let $p,q\in \P$. Suppose that $p$ is proper. Then
$q$ is  a dual of $p$ if and only if the mapping $r:{\mathcal K}_n\to \R^n$ defined by $r(K)= q(K^{p(K)})$ is itself an affine invariant point.
\end{lemma}
\par
\noindent
{\bf Proof.} By definition, if $q$ is the dual of $p$, then $r=p\in \P$.
Conversely, suppose that the map  $K\to q(K^{p(K)})$ is an affine invariant point. 
For all $K\in {\mathcal K}_n$, all linear, invertible maps  $T:\R^n\to \R^n$  and 
all $b\in \R^n$ we have by
 (\ref{adjoint})
$$
\Big(T(K)+b\Big)^{p\big(T(K)+b\big)} - p\Big(T(K)+b\Big) =\Big( T\big(K-p(K)\big)\Big)^{\circ}= T^{*-1} \Big(\big(K-p(K)\big)^{\circ}\Big).
$$
Therefore, and as $q$ and $r$ are affine invariant points, 
\begin{eqnarray*}
T^{*-1} \Big(q\big((K-p(K))^{\circ}\big)\Big) 
&= &q \Big(T^{*-1}\big((K-p(K))^{\circ}\big)\Big)\\
& =& q\left( \Big(T(K)+p(K)\Big)^{p\big(T(K)+p(K)\big)} - p\Big(T(K)+p(K)\Big) \right) \\
&=&
q\left( \Big(T(K)+p(K)\Big)^{p\big(T(K)+p(K)\big)} \right)- p\Big(T(K)+p(K)\Big) \\
&=&
r \Big(T(K)+p(K)\Big)- p\Big(T(K)+p(K)\Big)  
= T \Big(r(K)-  p(K)\Big) \\
&=& T\left(q\Big(\big(K-p(K)\big)^{\circ}\Big) \right).
\end{eqnarray*}
In particular, $T^{*-1} \Big(q\big((K-p(K))^{\circ}\big)\Big)=T\left(q\Big(\big(K-p(K)\big)^{\circ}\Big) \right)$ holds for $T= \lambda Id$ with $\lambda >1$. 
One has thus
$q\Big(\big(K-p(K)\big)^{\circ}\Big)=0$ for all $K$, and hence
$r(K)= p(K)$ for all $K\in \K$.\hskip 2mm$\square$
\vskip 3mm
We will show in Theorem \ref{theo:dual}  that if a proper affine invariant point $p$  has a dual, then this dual point is unique and proper.   We will then  show that $p$ is the unique dual   of $p^{\circ}$,  and hence  $p$ has a unique bidual point  which is $p$ itself.
\par
First,  we give a definition which will be useful to investigate duality.
\vskip 2mm
\begin{definition} \label{DEF}
Let $p\in \P$ be proper .
\par
\noindent
We say that  $p$ is injective if,  whenever $K_{1}, K_{2}\in \K$  satisfy
$K_{1}^{p(K_{1})}=K_{2}^{p(K_{2})}$, then $p(K_{1})=p(K_{2})$.
\par
\noindent
We say that $p$ is surjective if for  every
$C\in \K$, there exists $K\in \K$ such that
$C=K^{p(K)}$.
\par
\noindent
We say that $p$ is bijective if it is both injective and surjective.
\end{definition}

\vskip 2mm
The centroid, the Santal\'o point  the center of the John ellipsoid and the center of the L\"owner ellipsoid are examples of injective and surjective affine invariant points. More examples are given in 
Section \ref{example}.
\vskip 3mm

\noindent
\begin{remark}\label{RemInj} 
Let $p\in \P$ be proper and define $\phi_p:\K\to \K$ by 
\begin{equation} \label{phi}
\phi_p(K)=K^{p(K)}. 
\end{equation}
 It is easy to see that $\phi_p$ is  continuous. Moreover we have
\par
\noindent
(i)  {\it
 $p$ is injective, (surjective, bijective)  iff $\phi_p$ is injective,  (surjective, bijective)}.
\par
We address  the first statement. Let  $p$ be injective and suppose that $\phi_p(K_1)=\phi_p(K_2)$, i.e. $K_{1}^{p(K_{1})}=K_{2}^{p(K_{2})}$.
Then, by injectivity of $p$, $K_1=K_2$, i.e.  $\phi_p$ is injective.
Conversely,  let  $ \phi_p$ be  injective and suppose that $K_{1}^{p(K_{1})}=K_{2}^{p(K_{2})}$. The latter means exactly that $\phi_p(K_1)=\phi_p(K_2)$, 
and it follows from the injectivity of $\phi_p$, that $K_1=K_2$.
\vskip 1mm
\noindent
(ii)  {\it $q\in \P$ is dual of $p$ if and only if $q\circ\phi_p =p$}. 
\end{remark}
\vskip 3mm
The next  two lemmas characterize  injectivity and surjectivity.
\vskip 2mm

\begin{lemma}\label{surjective1}
Let $p$ be a proper affine invariant point. The following are equivalent.
\par
(i)  $p$ is surjective.
\par
(ii) For every $C$  in $\K$ there is a $z\in \inte(C)$ such that
$p((C-z)^{\circ}) =0$.
\end{lemma}
\par
\noindent
{\bf Proof.}  $p$ is surjective means that for all $C$  in $\K$ there is $K$  in $\K$  such that 
$$
C=K^{p(K)}  = \big(K-p(K)\big)^\circ +p(K), 
$$
or, equivalently,
$
\big(C-p(K)\big)^\circ = K- p(K).
$
This is equivalent to
$$
p\Big(\big(C-p(K)\big)^\circ\Big) = p(K)- p(K)=0.
$$
$\Box$
\vskip 3mm
\begin{lemma}\label{lemma:pinject}
Let  $p$ be a proper  affine invariant point.
Then the following are equivalent.
\newline
(i) $p$ is injective.
\newline
(ii) For all $C\in \K$, there exists at most one $z\in\inte(C)$ such that $p\big((C-z)^{\circ}\big)=0$.
\end{lemma}
\vskip 2mm
\noindent
{\bf Proof.}
$(i)\Rightarrow(ii)$  Suppose that there are $z_1$ and $z_2$ in $\inte(C)$ such that $p\big((C-z_1)^{\circ}\big)= p\big((C-z_2)^{\circ}\big)=0$. 
For $i = 1, 2$, we put $K_i=C^{z_i} = (C-z_i)^\circ +z_i$. Then $K_i - z_i = (C-z_i)^\circ$ and hence
$$
0 = p\big((C-z_i)^{\circ}\big) = p\big(K_i-z_i)\big) = p(K_i) -z_i, 
$$
and thus $p(K_i)= z_i$. By (\ref{bipolar}), $C=K_1^{z_1}= K_1^{p(K_1)}$ and $C=K_2^{z_2}= K_2^{p(K_2)}$. Injectivity of $p$ implies that $z_1=z_2$.
\par
\noindent
$(ii)\Rightarrow(i)$ 
Suppose $C=K_1^{p(K_1)}=K_2^{p(K_2)}$ for $K_1, K_2\in\K$. Then, for $i=1,2$,  $C-p(K_i)= \big((K_i-p(K_i)\big)^{\circ}$,  so that
$\big(C-p(K_i)\big)^{\circ} = K_i-p(K_i)$.  Hence  $C-p(K_i)$ has a bounded polar,  which means that $p(K_i)\in \inte(C)$. It follows that for $i=1, 2$, 
$$
p\Big(\big(C-p(K_i)\big)^{\circ}\Big) =p\Big( K_i-p(K_i)\Big)=0,
$$ 
and hence  by (ii) that $p(K_1)=p(K_2)$.  $\Box$
\vskip 3mm 
It will be useful to have  a new description of $(K^{\circ}-z)^{\circ}$ when  $0\in \inte(K)$ and $z\in \inte(K^{\circ})$.
Let thus $K\in \K$ be such that $0\in \inte(K)$. 
For  $z\in \inte(K^\circ)$,  we put 
\begin{equation}\label{K_z}
K_z=
\left\{ \frac{x}{1-\langle x, z\rangle} :  x\in K
\right\}.
\end{equation} 
In Lemma \ref{g0} we show that $K_{z}=(K^{\circ}-z)^{\circ}$.
\par
 It is easy to show
$$
|K_z|=\int_K \frac{dx}{(1-\langle x, z\rangle)^{n+1}} 
\hskip 10mm\mbox{and}\hskip 10mm
\lim_{z\to \partial K^{\circ}}|K_z|= +\infty.
$$ 
Moreover,  when $\lambda\to 1$, $|K_{\lambda z_0}|\to +\infty$ uniformly in $z_0\in \partial K^{\circ}$.
\par
For any ellipsoid ${\cal E}$ centered at $0$,   and every $z$ such that $h_{\cal E}(z)<1$, ${\cal E}_z$ is an ellipsoid.
For the Euclidean unit ball $B_2^n$, $(B_2^n)_z$ is an ellipsoid with center $\frac{z}{1-|z|^2}$ and
$$|(B_2^n)_z|= \frac{ |B_2^n|}{ (1-|z|^2)^{\frac{n+1}{2}}}.$$
\begin{remark} 
Therefore, a natural question to  ask is whether  ellipsoids are the unique bodies such that $K_z$ is centrally symmetric for any $z\in K^\circ$.
\end{remark}

\vskip 2mm
The next lemma relates $K_z$ to $(K^{\circ} -z)^{\circ}$.
\vskip 2mm

\begin{lemma}\label{g0} For all $K\in \K$ with $0\in \inte(K)$ and all $z\in\inte(K^{\circ})$, 
$$
(K^{\circ})^z-z= (K^{\circ} -z)^{\circ}= K_z.
$$
\end{lemma}
\vskip 2mm

By Lemma \ref{g0},
for all $z\in{\rm int}( K^\circ)$ and all $z'\in {\rm int}\big((K_z)^\circ\big)= \inte(K^{\circ}-z)$ 
$$
(K_z)_{z'}=\big( (K_z)^{\circ}-z'\big)^{\circ}= \big(K^{\circ}-(z+z')\big)^{\circ}=K_{z+z'}. 
$$
\noindent
{\bf Proof.}
The first equality follows from the definition. For the second one, observe that for $z\in\inte(K^{\circ})$, 
\begin{eqnarray*} 
(K^{\circ} -z)^{\circ}&= &\{x'\in \R^n:  \langle x',y-z\rangle \le 1\hbox{ for all } y\in K^{\circ}\}\\
&= & \{x'\in \R^n:  \langle x',y\rangle \le 1+\langle x',z\rangle \hbox{ for all } y\in K^{\circ}\}.
\end{eqnarray*} 
Since $0\in \inte\big(K^{\circ}\big)$,  such an $x'$  satisfies $1+\langle x',z\rangle > 0$, so that 
$$(K^{\circ} -z)^{\circ} =\left\{x'\in \R^n:  \left\langle \frac{x'} {1+\langle x',z\rangle}, y \right\rangle \le 1\hbox{ for all } y\in K^{\circ}\right\}.$$
Since $(K^{\circ})^{\circ}=K$, one has thus
$$(K^{\circ}-z)^{\circ}= \left\{x'\in \R^n:  \frac{x'} {1+\left\langle x',z \right\rangle}\in K\right\}. $$
Finally, observe that  
$x=\frac{x'} {1+\langle x',z\rangle}$  
if and only if  $x'=\frac{x} {1-\langle x,z\rangle}$.  \hskip 2mm  $\square $
\vskip 3mm 
Now we  show that  for a proper affine invariant point $p$ and for the centroid $g$, 
$p\big((K^{\circ}-\lambda x_0)^{\circ}\big)$ and $g\big((K^{\circ}-\lambda x_0)^{\circ}\big)$ 
have an analogous behavior when $\lambda\to 1$. 
To do so, we need several technical lemmas.
\vskip 2mm
\begin{lemma}\label{g1}
Let $K$ be a convex body containing $0$ in its interior.  Then there are constants $c>0$ and $0<\lambda_0<1$ such that for all 
$x_0^*\in \partial K^{\circ}$, 
all proper  $p\in \P$  and all  $\lambda_0<\lambda<1$, 
$$|p(K_{\lambda x_0^*})|\ge \left\langle p(K_{\lambda x_0^*}),\frac{x^*_0}{|x^*_0|} \right\rangle \ge \frac{c}{1-\lambda}.$$
\end{lemma}
\vskip 2mm
\noindent
{\bf Proof.} 
It is well known that for every $L\in {\mathcal K}_n$ one has that
$L-g(L)\subset n\big(g(L)-L\big), $ and  thus $L- (n+1)g(L)\subset-nL$. Hence,  for every $v\in S^{n-1}$, 
$$\langle g(L),v\rangle \ge \frac{1}{n+1}\left(h_L(v)-nh_L(-v)\right).$$
Since $p$ is proper, it follows from Proposition 1 of \cite{MeyerSchuettWerner2012} that for some $0\le \alpha<1$, one has
$p(L)-g(L)\in \alpha\big(L-g(L)\big)$.  Therefore,  for every $v\in S^{n-1}$, 
\begin{equation}\label{gleichung1}
\langle p(L),v\rangle  \ge (1-\alpha)\langle g(L),v\rangle -\alpha h_L(-v)\ge \frac{1-\alpha}{n+1}h_L(v)-\left(\frac{n(1-\alpha)}{n+1}+\alpha\right)h_L(-v).
\end{equation}
By Lemma \ref{g0}, for all $\lambda$ with $0<\lambda<1$
$$
K_{\lambda x_0^*} =\left\{ \frac{x}{ 1-\lambda \langle x_0^*, x\rangle} :  x\in K \right\}.
$$
Consequently, 
$$h_{K_{\lambda x_0^*} }\left(\frac{x^*_0}{|x^*_0|}\right)=
\frac{1}{|x^*_0|} \sup_{x\in K}  \frac{ \langle x^*_0,x\rangle } {1-\lambda \langle x,x_0^*\rangle}
= \frac{1}{|x^*_0|  (1-\lambda)  }$$
and
$$h_{K_{\lambda x_0^*} }\left(-\frac{x^*_0}{|x^*_0|} \right)=\frac{1}{\lambda |x^*_0|}
\sup_{x\in K} \frac{ -  \lambda \langle  x^*_0, x\rangle }{1-\lambda \langle x,x_0^*\rangle}
\le \frac{1}{\lambda|x^*_0|}.$$
Together with (\ref{gleichung1}) we get
$$\left\langle p(K_{\lambda x_0^*}),\frac{x^*_0}{|x^*_0|} \right\rangle \ge \frac{1}{|x^*_0|}\left(\frac{1-\alpha}{(1-\lambda)(n+1)}
                                     - \frac{1}{\lambda}\left(\frac{n(1-\alpha)}{n+1}+\alpha\right)\right).$$
We choose $$c= \frac{ 1-\alpha}{ 2(n+1)  \  \max_{x_0^*\in \partial K^{\circ}} |x_0^*|}  $$  and we obtain the result, choosing a 
big enough $\lambda_0$, for all $\lambda$ with $0<\lambda_0 \leq \lambda<1$.  $\Box$
\vskip 3mm
\begin{lemma}\label{g12}
 For all $K\in\K$ with a $C_2^+$-boundary there exists $c'>0$ such that for all $(x_0^*,x_0)\in \partial K^{\circ}\times \partial K$  satisfying
$\langle x_0^*, x_0\rangle =1$ and all $x\in K$
 $$|x-x_0|^2\le 
c'(1-\langle x_0^*,x\rangle).$$
\end{lemma}
\vskip 2mm
\noindent
{\bf Proof.} By the $C_2^+$ hypothesis, there exists an $R>0$ such that   for all $(x_0^*,x_0)\in \partial K^{\circ}\times \partial K$ with $\langle x_0^*,x_0\rangle=1$,  
$$K \subset x_0-R\frac{x_0^*}{|x_0^*|} +RB_2^n.$$
Hence, for   all $x\in K$, $|x_0-x -R\frac{x_0^*}{|x_0^*|}|^2\le R^2$, or, 
$$|x_0-x|^2\le \frac{2R}{|x_0^*|}\langle x_0^*, x_0-x\rangle= \frac{2R}{|x_0^*|} (1-\langle x_0^*, x\rangle).$$ We take  $c'= \frac{2R}{\min_{y^*\in \partial K^{\circ}}|y^*|}$.\hskip 2mm $\square$
\vskip 3mm
\begin{lemma}\label{g123}
For all $K\in\K$ with a $C_2^+$-boundary and $0 \in \inte(K)$ there exists $d>0$ such that for all $(x_0^*,x_0)\in \partial K^{\circ}\times \partial K$ with 
$\langle x_0^*, x_0\rangle =1$, for all $v\in S^{n-1}$ with 
$\langle v,x_0\rangle=0$, for all $0 < \lambda < 1$ and for all $y\in K_{\lambda x_0^*}$, 
 $$
 \langle v,y\rangle\le \frac{d}{2\sqrt{\lambda(1-\lambda)}}.
 $$
\end{lemma}
\vskip 2mm\noindent
{\bf Proof.} 
Let  $0 < \lambda <1$ and let $y\in K_{\lambda x_0^*}$. Then 
$y=\frac{x}{1-\lambda \langle x_0^*,x\rangle}$  for some $x\in K$. By Lemma \ref{g12}
and $\langle v,x_{0}\rangle=0$
$$
|\langle v,y\rangle|= \frac{|\langle x,v\rangle|}{1-\lambda \langle x_0^*,x\rangle}\le \frac{|x_ 0-x|}{1-\lambda \langle x_0^*,x\rangle}
\le\frac{   \sqrt{ c (1-\langle x_0^*,x\rangle) } }  {1-\lambda \langle x_0^*,x\rangle }\le \frac{\sqrt{c}}{2\sqrt{\lambda(1-\lambda)}}.
$$
The last inequality follows as for $t<1$ and $0<\lambda<1$
 $$
 \frac{\sqrt{1-t}}{1-\lambda t}\le \frac{1}{2\sqrt{\lambda(1-\lambda)}}.
 $$
 $\Box$
\vskip 3mm
\begin{proposition}\label{g1234}
 Let $K\in \K$ be in $C_2^+$ and such that $0 \in \inte(K)$.  Let $p$ be a proper affine invariant point. Then there is a constant $C$ such that for all $\lambda$ with 
 $\frac{1}{2}\leq\lambda<1$, for all $x^*_0\in \partial K^{\circ}$ and all $v\in S^{n-1}$
 with $\langle x_{0},v\rangle=0$
 \begin{equation}\label{hilfslimes}
\frac{ \langle p_{\lambda},  \frac{x_0}{|x_0|}\rangle }
 { |\langle p_{\lambda},v\rangle|}
\ge\frac{C}{\sqrt{1-\lambda}},
\end{equation}
where $x_{0}\in\partial K$ is the unique point with $\langle x_{0},x_{0}^{*}\rangle=1$.
In particular, uniformly in $x^*_0\in \partial K^{\circ}$, 
$$
\lim_{\lambda \to 1} \ \frac{p\big((K^{\circ}-\lambda x_0^*)^{\circ}\big)}
{\big|p\big((K^{\circ}-\lambda x_0^*)^{\circ}\big)\big|}= \frac{x_0}{|x_0|}.
$$
\end{proposition}
\vskip 2mm
\noindent
{\bf Proof.}
Let $u= \frac{x_0}{|x_0|}\in S^{n-1}$, $u^*=\frac{x^*_0}{|x^*_0|} \in S^{n-1}$ and for $0 < \lambda <1$ let $p_{\lambda} =p\big((K^{\circ}-\lambda x_0^*)^{\circ}\big)$.  To show (\ref{hilfslimes}), we choose
$w\in S^{n-1}$ such that
$$
u^*=\langle u^*,u\rangle u +\sqrt{1-\langle u^*,u\rangle^2} w.
$$ 
Please note that $\langle w,u\rangle =0$. Then
$$ \langle p_{\lambda}, u\rangle =\frac{1}{\langle u^*,u\rangle} \left(\langle p_{\lambda}, u^*\rangle -\sqrt{1-\langle u^*,u\rangle^2} \langle p_{\lambda},w\rangle\right).
$$
By Lemmas \ref{g1} and \ref{g123}, for $\lambda $ big enough, $p$ is proper
and all $v\in S^{n-1}$ with $\langle x_{0},v\rangle=0$
$$
\frac{ \langle p_{\lambda}, u\rangle } { |\langle p_{\lambda},v\rangle|}
\ge \frac{  \frac{1}{\langle u^*,u\rangle}\big(     \frac{c}{1-\lambda}- \sqrt{1-\langle u^*,u\rangle^2}\frac{d}{2\sqrt{\lambda(1-\lambda)}}\big)}
{\frac{d}{2\sqrt{\lambda(1-\lambda)}} },
$$
where $d$ denotes the constant from Lemma \ref{g123}.
Since 
$$\frac{1}{\langle u^*,u\rangle }\ge 
\min_{x_0^*\in  \partial K^{\circ}}{|x^*_0|\ |x_0|} >0
$$
there is a constant $C$ such that for all $\lambda$ with $\frac{1}{2}\leq\lambda<1$
and all $x_{0}\in\partial K$ and $x_{0}^{*}\in\partial K^{\circ}$ 
with $\langle x_{0},x_{0}^{*}\rangle=1$
$$
\frac{ \langle p_{\lambda}, u\rangle } { |\langle p_{\lambda},v\rangle|}
\ge\frac{C}{\sqrt{1-\lambda}}.
$$
By  Lemma \ref{g1}, $|p_{\lambda} |\to +\infty$ when $\lambda\to 1$.
We write $p_\lambda= \langle p_{\lambda}, u\rangle  u + \left(p_\lambda -  \langle p_{\lambda}, u\rangle u  \right)$ and put $v= \frac{p_\lambda -  \langle p_{\lambda}, u\rangle u}{|p_\lambda +  \langle p_{\lambda}, u\rangle u|}$. Then $v$ is orthogonal to $u$ and
\begin{eqnarray*}
\frac{p_\lambda}{|p_\lambda|} 
&=& \frac{\langle p_{\lambda}, u\rangle u 
+ \langle p_{\lambda}, v\rangle v}{\left[(\langle p_{\lambda}, u\rangle)^2  
+ (\langle p_{\lambda}, v\rangle)^2 \right]^\frac{1}{2}}
 = \frac{u}{
\left[1  + \frac{(\langle p_{\lambda}, v\rangle)^2}{(\langle p_{\lambda}, u\rangle)^2} \right]^\frac{1}{2}} + \frac{v}{\left[1  + \frac{(\langle p_{\lambda}, u\rangle)^2}{(\langle p_{\lambda}, v\rangle)^2} 
\right]^\frac{1}{2}},
\end{eqnarray*}
which converges to $u$ if (\ref{hilfslimes}) holds. $\Box$
\vskip 3mm
\begin{proposition}\label{g6}
Let $p\in \P$ be proper. For   $C\in\K$,  let $F:\inte(C)\to \R^n$ be  the mapping defined by 
$$F(z)=p\big((C-z)^{\circ}\big).$$ 
Then  $F$ is surjective.
\end{proposition}
\vskip 2mm \noindent
{\bf Proof.}
We fix  $z_0\in\inte(C)$ and put $L=(C-z_0)^{\circ}$. Then   $0\in\inte(L)$ and
$L^{\circ}=C-z_0$, so that the statement of the proposition is equivalent to:
\par
For any $K\in \K$ such that 
$0\in\inte(K)$,  the mapping $F:\inte(K^{\circ})\to \R^n$ defined by
$F(y^*)=p\big((K^{\circ}-y^*)^{\circ}\big)$
is surjective.
\par
We shall prove it in this form.
\par
We first treat the case when the body $K$ is $C_2^+$.  Suppose that $F$ is not surjective. Then for some $x\in \R^n$, one has $F(y^*)\not=x$ for every $y^*\in \inte(K^{\circ})$. 
We define  the function $G: \inte(K^{\circ})\to S^{n-1}$ by
$$G(y^*)=\frac{ F(y^*)-x} {|F(y^*)-x|}.$$
By Lemma \ref{g1},  $|F(y^*)|\to \infty$ for $y^*\to \partial K^{\circ}$.  By this and  by Proposition \ref{g1234},  we have for all $y^*\in \inte(K^{\circ})$,
$$
\lim_{ y^*\to\partial K^{\circ}} \frac{ F(y^*)-x} {|F(y^*)-x|} = \frac{y}{|y|} ,
$$
 where $y$ is the unique point in $\partial K$ such that $\langle y^*,y\rangle=1$.
Thus $G$ can be extended  to a continuous function $H:  K^{\circ}\to S^{n-1}$, setting 
\begin{eqnarray*}
H(y^*) = \left\{
\begin{array}{cc}  G(y^*) &  \hskip 4mm  \text{if } \   y^*\in \inte(K^{\circ}) \\
 \frac{y}{|y|}   &    \text{if }   \  y^* \in \partial K^{\circ}, 
  \end{array}
 \right.
\end{eqnarray*}
Indeed, $H$ is continuous on $ \partial K^\circ$ since $\partial K$ is $C^2_+$.
We also define  a continuous function $\theta: B_2^n\to K^{\circ}$ by 
\begin{eqnarray*}
\theta(z)= \left\{
\begin{array}{cc} 
0&  \hskip 4mm  \text{if } \  z=0\\
 \frac{ |z|}{h_K(z)} z  &    \text{if }   \  z\in B_2^n\setminus \{0\}. 
   \end{array}
 \right.
\end{eqnarray*}
Then $\phi=H\circ \theta: B_2^n\to S^{n-1}$ is continuous.
It follows that the function $\psi:B_2^n\to B_2^n$ defined by
$$\psi(z)=\frac{ z-\phi(z)}{2}$$
is also continuous. By the Brouwer fixed point theorem (see e.g. \cite{Rudin}),  for some  $z_0\in B_2^n$, one has $\psi(z_0)=z_0$, so that
 $$
 z_0=-\phi(z_0)\in S^{n-1}.
 $$ 
 $z_{0}\ne0$ since $z_{0}\in S^{n-1}$. Therefore
 $$
 \theta(z_{0})=\frac{|z_0| \  z_0}{h_K(z_0)}
 $$ 
and consequently $x_0^*=\theta(z_0)=\frac{|z_0| \  z_0}{h_K(z_0)}\in \partial K^{\circ}$. 
One gets
$$
\phi(z_0)=H(x_0^*)=\frac{x_0}{|x_0|} ,
$$
where $x_0$ is the unique point in $\partial K$ such that $\langle x_0^*,x_0\rangle =1$. Since
$$ z_0=\frac {x_0^*}{| x_0^*|}=-\phi(z_0)= -\frac{x_0}{|x_0|}, $$
we get
$$\langle z_0,z_0\rangle=\langle z_0,-\phi(z_0)\rangle=-\left\langle \frac{x_0}{|x_0|},\frac {x_0^*}{| x_0^*|}\right\rangle
= -\frac{1}{|x_0| \ |x_0^*|} < 0,$$
which is absurd.
\par
Now we treat the general case. Let $K$ be a convex body such that $0\in \inte(K)$. Then, there exists a sequence $(K_m)$ of $C_2^+$ convex bodies converging to $K$ in the 
Hausdorff metric. For $m$ big enough, one has $0\in \inte(K_m)$. Fix $x\in\R^n$. By above, there exists $y_m^*\in \inte(K_m^{\circ})$ such
that $p\big((K_m^{\circ}-y_m^*)^{\circ}\big)=x$. Since $K_m\to K$, $K_m^{\circ} \to K^{\circ}$. Let 
$y_{m_k}^*\to y^*$ be a converging subsequence of the sequence $(y_m^*)_{m \in \mathbb{N}}$. Then
it is easy to see that $y^*\in\inte(K^{\circ})$, and,  by continuity,
$p\big((K^{\circ}-y^*)^{\circ}\big)=x$.\hskip  2mm $\square$

\vskip 3mm
\section{The main theorems.}\label{main}

The first  theorem in this section  follows immediately from Lemma
\ref{surjective1} and Proposition \ref{g6}. Observe however that the statement of Proposition  \ref{g6} is  stronger than surjectivity of 
all proper $p\in \P$.
\vskip 2mm
\begin{theorem}\label{theo:surjective}
Every proper affine invariant point  is surjective.
\end{theorem}
\vskip 1mm
While every proper affine invariant point is surjective, for injectivity this is not the case. This is the content of the next theorem.
\vskip 2mm 

\begin{theorem}\label{theo:nichtinjective}
For $n \geq 2$, there exists  $p  \in \mathfrak{P}_n$ such that $p$ is not injective.
\end{theorem}

\vskip 1mm
We postpone the proof of Theorem \ref{theo:nichtinjective} to the end of this section. First we apply Theorem \ref{theo:surjective}  to obtain the following result.
\vskip 3mm

\begin{theorem} \label{theo:dual}
 Let $p\in \P$ be  proper. Then the following assertions are equivalent.
 \par
 (i) $p$ has a dual $q$.
\par 
 (ii) There exists a proper $r\in \P$ such that $p$ is a dual of $r$.
\par 
 (iii) $p$ is injective.
 \par
 \noindent
Moreover, if one of these assertions hold, then $p$ has a unique dual point $p^{\circ}$.  $p^{\circ}$ is proper,  $r=p^{\circ}$
and $p^{\circ}$ has a unique dual which is $p$. 
\end{theorem}
\vskip 2mm
\noindent
\subsection{ Proof of Theorem \ref{theo:dual}.}
We shall need more lemmas.
\vskip 2mm
\begin{lemma} \label{useful} Let $K\in \K$ and suppose that  $K-y\subset \beta(y-K)$ for some $\beta> 1$ and some $y\in \R^n$.  
For all  $x\in \R^n$ and  for all real numbers $\alpha$ and $\gamma$ such that  $0<\alpha <1<\gamma$ the following assertions hold.
\par
(i)  If  $x-y\in \alpha(K-y)$, then $K-x\subset \frac{\beta+\alpha}{1-\alpha} \  (x-K)$.
\par
(ii)  If $K-x \subset \gamma (x-K)$, then $x-y\in \frac{\beta\gamma-1}{\beta(\gamma +1)} \  (K-y)$.
\end{lemma}
\vskip 2mm
\noindent
{\bf Proof.} We may  assume that $y=0$. Otherwise we consider the body $K'=K-y$.
\par
\noindent
(i)  Suppose $x\in \alpha K$. Let $\delta= \frac{\beta +\alpha}{1-\alpha}$. We need to prove that 
$$(\delta +1)x-K\subset \delta K.$$
With the assumption in  (i), the inclusion $-K\subset \beta K$ and the convexity of $K$, we get
$$(\delta+1)x-K\subset \alpha(\delta +1)K +\beta K=\left(\alpha\left(1+\frac{\beta+\alpha}{1-\alpha} \right)+\beta\right)K=\delta K.$$
\par
\noindent
(ii)  Suppose that $K-x \subset  \gamma (x-K)$. Then $(\gamma +1)x -K\subset \gamma K.$
Now we use that  $\frac{K} {\beta}\subset -K$, divide by $\gamma$ and 
get with $y=\frac{\gamma +1}{\gamma} x$ and $t=\frac{1}{\beta \gamma}$ that 
$$y +tK\subset K.$$
Since $K$ is bounded and closed, it follows that $y\in (1-t)K$,
so that 
$$x=\frac{\gamma}{\gamma +1} y\in \frac{\beta \gamma- 1}{\beta(\gamma +1)}K.\hskip 2mm \square$$

\vskip 3mm \noindent
\begin{lemma}\label{veryuseful}
Let $L\in \K$ and let  $0<r \leq R < \infty$   be such that $rB_2^n\subset L\subset R B_2^n$. Suppose that there is $x\in \R^n$  and $\gamma \ge 1$ such that
$L-x\subset \gamma (x-L)$. Then
$$\frac{2 r}{\gamma +1}B_2^n\subset L-x \subset \frac{2\gamma R}{\gamma +1}B_2^n$$
\end{lemma}
\vskip 0mm
\noindent
{\bf Proof.} 
As $L-x\subset \gamma (x-L)$, $$2rB_2^n \subset L-L= L-x+x-L\subset \gamma (x-L)+x-L=(\gamma +1)(x-L).$$
This gives the first inclusion. The second one is obtained from
$$
(\gamma +1)(L-x)=\gamma(L-x)+L-x\subset \gamma(L-x)-\gamma (L-x)= \gamma(L-L)\subset 2\gamma R B_2^n.
$$
$\Box$

\vskip 3mm
For the next lemma, recall (\ref{phi}) where  we introduced the notation   $\phi_p (K)= K^{p(K)}$.
\vskip 2mm

\begin{lemma}\label{veryveryuseful} 
Let $p$ be a proper affine invariant  point and suppose that the mapping $\phi_p$
is bijective. Let $q:\K\to \R^n$ be defined
by $q(L)= p(K)$ whenever $L=K^{p(K)}$.
Then $\phi_p$ is a homeomorphism and the mapping  $q$ is continuous.
\end{lemma}
\vskip 2mm

\noindent
{\bf Proof.} It is clear that $\phi_p$ is continuous. We want to show that $\phi_{p}^{-1}$
is continuous.
We use now that the inverse of a continuous, bijective map
between locally compact Hausdorff spaces is continuous if the inverse image of any compact set
is compact.
Since $\K$ endowed with the Hausdorff metric is locally compact, it is enough to verify that for any compact subset ${\cal L}$ of $\K$, $(\phi_p)^{-1} ({\cal L})$ is compact.
\par
By Proposition 1 of \cite{MeyerSchuettWerner2012}, one has for every $K\in \K$ that
$$p(K)-g(K)\subset \alpha \big(( K-g(K)\big), $$ 
for some $0< \alpha <1$. Here, $g$ denotes the centroid.
It is well known that 
$$ K-g(K)\subset n\big(g(K)-K\big).$$
It follows from Lemma \ref{useful} (i) 
$$K-p(K)\subset  \frac{n+\alpha}{1-\alpha} \big(p(K)-K\big),$$ whence $$\big(K-p(K)\big)^{\circ}\subset -  \frac{n+\alpha}{1-\alpha} \big(K-p(K)\big)^{\circ}.$$ 
Now
$\phi_p(K)= \big((K-p(K)\big)^{\circ} +p(K)$, so that 
\begin{equation} \label{Nummer11}
\phi_p(K)-p(K)\subset  \frac{n+\alpha}{1-\alpha} \big(p(K)-\phi_p(K)\big).
\end{equation}
Let  ${\cal L}$ be a compact subset of $\K$. By affine invariance, we may suppose that there are  $0<r\leq R < \infty$ such that for every $L\in {\cal L}$
$$rB_2^n\subset L\subset RB_2^n.$$
Let $L\in {\cal L}$ and $K=(\phi_p)^{-1}(L)$,  that is,  
$$L-p(K)=\phi_p(K)-p(K)
\subset \frac{n+\alpha}{1-\alpha}  \left(p(K) - \phi_p(K) \right) = \frac{n+\alpha}{1-\alpha}\left( p(K) - L\right) .$$ 
The last inclusion follows from  (\ref{Nummer11}).
Therefore, by Lemma \ref{veryuseful}  with $\gamma= \frac{n+\alpha}{1-\alpha}$, 
$$\frac{2 r}{\gamma +1}B_2^n\subset L-p(K) \subset \frac{2\gamma R}{\gamma +1}B_2^n.$$
By duality $$\frac{\gamma +1}{2\gamma R}B_2^n\subset K-p(K) \subset \frac{\gamma +1}{2 r}B_2^n.$$
As $p$ is proper,  $p(K)\in L\subset RB_2^n$. Hence  we get that for some $0<c<d$, one has 
$$
(\phi_p)^{-1}({\cal L})\subset {\cal L}'=\{K\in\K: \ x+cB_2^n \subset K\subset d B_2^n\hbox{ for some }x\in \R^n \}.
$$
Now, it is easily seen that
${\cal L}'$ is a compact subset of $\K$, and consequently $(\phi_p)^{-1}({\cal L})$ is as 
a closed subset of a compact set compact.
$\Box$

\vskip 2mm\noindent
{\bf Proof of Theorem \ref{theo:dual}.}

(i)$\implies$(iii) Suppose that  $K_1, K_2\in \K$ satisfy $K_1^{p(K_1)} =K_2^{p(K_2)}$. Then, by the definition of  dual point, 
$$
p(K_1) =  q\left(K_1^{p(K_1)}\right) =q\left(K_1^{p(K_1)}\right) =p(K_2).
$$
\par
(iii) $\implies$(i) By Theorem \ref{theo:surjective}, $p$ is surjective and consequently bijective. 
By  Remark  \ref{RemInj}(i), this is equivalent that $\phi_p$ is bijective.   Lemma \ref{veryveryuseful} then gives
$$
q(L) = q(K^{p(K)}) = p(K)
$$
and thus $q$ is the dual of $p$.
\par
(i) $\implies$ (ii) We first show that  $q$ is proper. By Theorem \ref{theo:surjective},  $p$ is surjective. Hence there exists $K\in \K$ if $C\in \K$ such that  $C=K^{p(K)}$. 
As $q$ is the dual of $p$,  one has $q( C )= q(K^{p(K)}) = p(K)$ and  thus $C=K^{p(K)}= K^{q( C )}$.  By the bipolar theorem  (\ref{bipolar}),  $C^{q( C )} = \left(K^{q( C )}\right)^{q( C )}=K$,  which proves that $q(C)\in \inte( C )$.
Next we show:
\begin{equation}\label{dual-hinundher}
\text{ If a proper } p \in \mathcal{P}_n  \text{  has a dual } q, \text{ then } p \text{ is a dual of} 
\hskip 1mm q. 
\end{equation}
Again, by Theorem \ref{theo:surjective},  for all  $L\in \K$, one has $L=K^{p(K)}$ for some $K\in \K$. Thus, using the definition of duality,   $p(K)= q\left(K^{p(K)}\right)$ and the bipolar theorem (\ref{bipolar}),
$$
L^{q(L)} = \left(K^{p(K)}\right)^{q\left(K^{p(K)}\right)}  = \left( K^{p(K)}\right) ^{p(K)} = K, 
$$
and hence  $p\left( L^{q(L)}\right) = p(K) = q(L)$, which proves that $p$ is a dual of $q$. 
Thus we can take $r=q$ in (ii). 
\par
(ii) $\implies$ (i) This follows immediately with (\ref{dual-hinundher}).  
$\Box$
\vskip 3mm
\subsection{A product  mapping on affine invariant points.}
The next definition  will turn out to be useful to characterize duality. There, $I_n:\P\to \P$ denotes the identity map. Recall also  the mapping $\phi_p: \K\to \K$   defined in (\ref{phi}) by
$\phi_p(K)= K^{p(K)}$. 

\vskip 2mm 

\begin{definition} \label{[p,q]}
Let  $(p,q)\in \P\times \P$ be  such $p$ and $q$ are proper.  We define  the  map $[p,q]:\P\to \P$ by
$$[p,q](r)= r\circ\phi_q\circ\phi_p-q\circ \phi_p+p.
$$
\end{definition}
\vskip 2mm 
We show first that   $[p,q](r)$ is indeed in $ \P$.
\par
\noindent
Since $\phi_p$ is continuous for all proper $p\in \P$,
$[p,q](r)$ is continuous for all $r\in \P$. 
\par
\noindent
Now we address affine invariance.
 Let $b\in \R^n$,  $T:\R^n \to \R^n$  be a bijective linear map 
and $S:\R^n\to \R^n$ be defined by $Sx=Tx+b$.  
Put $L=\big(K-p(K)\big)^{\circ}$. By definition,  for all $K\in \K$,
$$[p,q](r)(K)= r\Big(\big(L-q(L)\big)^{\circ}\Big) +p(K).
$$
We put $L'= \big( S(K)-p(S(K))\big)^{\circ}$. Then with (\ref{adjoint}),  
$$L'= \big( T(K)-p(T(K))\big)^{\circ}= \Big (T\big(K-p(K)\big)\Big)^{\circ} = \left(T\left(L^\circ\right) \right)^\circ =  T^{*-1}(L).
$$
Thus, again with (\ref{adjoint}),  
$$\big(L'-q(L')\big)^{\circ}= \bigg(T^{*-1}(L) - q \left(T^{*-1}(L)\right) \bigg)^\circ=\bigg(T^{*-1}\Big(L-q(L)\Big)\bigg)^\circ
=T\Big( \big(L-q(L)\big)^{\circ}\bigg), $$
so that 
\begin{eqnarray*} 
[p,q](r)(S(K))&=& r\Big(\big(L'-q(L')\big)^{\circ} \Big)+p(S(K))= r\bigg(T\Big( \big(L-q(L)\big)^{\circ}\Big)\bigg)+p(T(K))+b\\
&=& T\bigg( r\Big(\big(L-q(L)\big)^{\circ}\Big) +p(K)\bigg)+b= S\big([p,q](r)(K)\big).
\end{eqnarray*} 
\vskip 3mm 
With the aid of $[p,q]$, we give another characterization of dual affine invariant  points.
\vskip 2mm

\begin{proposition}\label{[p,q]=I}
Let $p,q\in \P$ be proper.  The following are equivalent.
\par
(i) $[p,q]=I_n$. 
\par
(ii) $[p,q](p)=p$.
\par
(iii) $q$ is the dual of $p$.
\par
(iv) $p$ is the dual of $q$.
\end{proposition}
\noindent
{\bf Proof.}
It is clear that (i)$\implies$(ii) and by Theorem \ref{theo:dual},  (iii) is equivalent to (iv).
We now show the remaining implications.
\par
(ii) $\implies$(iv)  By Theorem \ref{theo:surjective}, $p$ is surjective. Hence  for every $L\in \K$ there is  $K\in \K$ such that $L=\phi_p(K)$.  By (ii), 
\begin{equation*}
p(K)= [p,q] ( p )  (K) = (p\circ\phi_q\circ\phi_p)(K) - (q\circ \phi_p)(K)+p(K), 
\end{equation*}
which is equivalent to 
\begin{equation*}
p\left(\left(\phi_p(K)\right)^{q\left(\phi_p(K)\right)}\right) - q\left(\phi_p(K)\right) =0, 
\end{equation*}
and again,  equivalent to 
\begin{equation*}
p\left(L^{q\left(L\right)}\right) = q\left(L\right) .
\end{equation*}
By the definition, this means that $p$ is the dual of $q$.  
\par
(iii)$\implies$(i) If $q$ is the dual of $p$, then $q(K^{p(K)})=p(K)$ for every $K\in \K$.
This  means that $q\circ\phi_p=p$.  It follows that
$$
(\phi_q\circ\phi_p)(K)= (\phi_p(K)\big)^{q\big(\phi_p(K)\big)}=(\phi_p(K)\big)^{p(K)}
=\left(K^{p(K)}\right)^{p(K)}= K.
$$
The last equality follows from the Bipolar Theorem.
One has thus
$$
[p,q](r)=r\circ\phi_q\circ\phi_p -q\circ\phi_p+p= r-p+p=r.
$$
$\Box$
\vskip 3mm 
The next proposition describes the product $[p_1,q_1]\circ[p_2,q_2]$ in some special cases.
Note also  that 
\begin{equation}\label{Formel}
\phi_p\circ\phi_{p^{\circ}}= \phi_{p^{\circ}}\circ\phi_p = I_n, \  \  
 p=p^{\circ}\circ \phi_p \  \  {\text and} \  \  p^{\circ} =p\circ\phi_{p^{\circ}}.
 \end{equation}
\vskip 2mm \noindent
\begin{proposition}\label{composition}
Let $p,r,s\in\P$ be proper and suppose that $p$ has a dual $p^{\circ}$. Then
$$
[r,p^{\circ}]\circ[p,s]=[r,s].
$$
In particular,
if $p,q\in \P$ are proper and have dual points $p^{\circ}$ and $q^{\circ}$, then
$$
[q^{\circ},p^{\circ}]\circ [p,q]= I_n.
$$
\end{proposition}
\vskip 2mm
\noindent
{\bf Proof.} 
Let $t\in \P$. Then with  (\ref{Formel}), 
\begin{eqnarray*} 
 && ([r,p^{\circ}]\circ [p,s])(t)   \\
&&= [p,s](t)\circ\phi_{p^{\circ}}  \circ\phi_{r}- p^{\circ}\circ\phi_r+r\\
&&=  (t\circ \phi_s\circ\phi_p-s\circ\phi_p +p)  
\circ\phi_{p^{\circ}}  \circ\phi_{r}- p^{\circ}\circ\phi_r+r\\
&&= t\circ \phi_s\circ\phi_p\circ\phi_{p^{\circ}}  \circ\phi_{r}-s\circ\phi_p\circ\phi_{p^{\circ}}  \circ\phi_{r}+p\circ\phi_{p^{\circ}}  \circ\phi_{r}
-p^{\circ}\circ\phi_r+r\\
&&= t\circ \phi_s  \circ\phi_{r}-s\circ\phi_r+p^{\circ}\circ\phi_{r}-p^{\circ}\circ\phi_r+r\\
&&= t\circ \phi_s  \circ\phi_{r}-s\circ\phi_r+r= [r,s](t).
\end{eqnarray*} 
Therefore $[r,p^{\circ}]\circ [p,s]=[r,s]$ and
$[q^{\circ},p^{\circ}]\circ [p,q]= [q^{\circ},q]=I_n$. The last equality follows by
Proposition \ref{[p,q]=I}.
$\Box$
\vskip 2mm

 \noindent
\begin{remark} It would be interesting to know more about $[p,p]$ and  
 $$[p,p]^k(p):=[p,p]\circ\dots\circ [p,p]) (p) \hbox{ , $k$ times, $k\ge 1$.}$$
Is there a limit for $k\to \infty$ ?
\end{remark}
\vskip 3mm

Let $\A$ be as in Definition \ref{aism} and let  $A\in \A$. Let $p$, $a$ and $b$ be in $\P$. Suppose in addition  that $p$ is proper  and  that 
$a(M)\in \inte(A(M))$ for any $M\in \K$.
For $K\in \K$,  we define 
\begin{equation} \label{BK}
B(K) = \bigg(A\Big(\big(K-p(K)\big)^{\circ} \Big)-a\Big(\big(K-p(K)\big)^{\circ} \Big) \bigg)^{\circ}+ b(K).
\end{equation}
\par
Then $B$ is an affine invariant set mapping, i.e.  $B\in \A$. We now show this.
\par
\noindent
It is clear from the hypothesis on $a$ that 
$B(K)\in \K$ and that $B$ is continuous.  
We prove next  that $B$ is an affine invariant mapping.
To do so, fix $K\in \K$ and put $M=\big(K-p(K)\big)^{\circ}$. Then
\begin{equation} \label{beschraenkt}
\big(B(K)-b(K)\big)^{\circ}= A(M)-a(M),
\end{equation}
which shows that $b(K) \in \inte\left(B(K)\right)$.
Let $T:\R^n\to \R^n$ be linear and one-to-one, $c\in \R^n$ and $S=T+c$. With (\ref{adjoint}), 
\begin{eqnarray*} \label{Tstern}
\Big(S(K) - p\big(S(K)\big)\Big)^{\circ} &=& \Big(T(K)+c-p\big(T(K)+c\big)\Big)^{\circ}
= \Big(T\big(K-p(K)\big)\Big)^{\circ}\\&=&T^{*-1}\big ((K-p(K))^\circ\big) =
T^{*-1}(M).
\end{eqnarray*}
We get from this and (\ref{beschraenkt}) with $S(K)$ instead of $K$, 
\begin{eqnarray*}
\Big(B\big(S(K)\big)-b\big(S(K)\big)\Big)^{\circ}
= A\big(T^{*-1}(M)\big)-a\big(T^{*-1}(M)\big) 
=T^{*-1}\big(A(M)-a(M)\big).
\end{eqnarray*}
It follows with (\ref{adjoint}) and (\ref{beschraenkt}) that 
\begin{eqnarray*}
B\big(S(K)\big)-b\big(S(K)\big)
&=& T\Big(\big(A(M)-a(M)\big)^{\circ} \Big)
= T\big(B(K)-b(K)\big) \\
&=& T\big(B(K)\big) +c - T\big(b(K)\big) -c = S\big(B(K)\big)-b\big(S(K)\big), 
\end{eqnarray*}
so that 
$B\big(S(K)\big) = S\big(B(K)\big)$.  \vskip 2mm
Suppose now that $p$ has a dual $p^{\circ}$. By Theorem \ref{theo:dual},  $p^\circ$ is surjective.  Hence for all $K \in \K$ there is 
$L\in \K$ such that $ K=L^{p^{\circ}(L)} $. Using also  (\ref{Formel}), this implies that 
$$
p(K)= p\left( L^{p^\circ(L)}\right) = \left(p \circ \phi_{p^\circ}\right)(L) = p^\circ(L), 
$$
and
$$
\big(K-p(K)\big)^{\circ}=  \phi_p(K) -p(K) =  \left(\phi_p \circ \phi_{p^\circ} \right)(L)= L-p^{\circ}(L).
$$
One has then, also with (\ref{BK}), 
\begin{eqnarray*}
B\big(K\big)-b(K) &=& \bigg(A\Big(\big(K-p(K)\big)^{\circ} \Big)-a\Big(\big(K-p(K)\big)^{\circ} \Big) \bigg)^{\circ} \\
&=&  \big((A(L)-a(L)\big)^{\circ}.
\end{eqnarray*}
It follows that 
$$A(L)= \Big(B\big(K\big)-b(K)\Big)^{\circ}+a(L)= \bigg(  B\Big(\big(L-p^{\circ}(L)\big)^{\circ}\Big)-b\Big(\big(L-p^{\circ}(L)\big)^{\circ}\Big)\bigg)^{\circ}+a(L).$$

\vskip 3mm 
Thus  we have proved  the following proposition.
\vskip 2mm
\begin{proposition}
Let $A\in \A$ and let  $p$, $a$ and $b$ be in $\P$. Suppose in addition  that $p$ is proper  and  that 
$a(M)\in \inte(A(M))$ for any $M\in \K$. 
\newline
For $K\in \K$, $B$ defined by
$$B(K) = \bigg(A\Big(\big(K-p(K)\big)^{\circ} \Big)-a\Big(\big(K-p(K)\big)^{\circ} \Big) \bigg)^{\circ}+ b(K)$$
is an affine invariant set mapping.
\par
If  $p$ has a dual point $p^{\circ}$, then $b(M)\in \inte\big(B(M)\big)$ for any $M\in \K$ and   for any $L\in \K$,
$$A(L)=  \bigg(  B\Big(\big(L-p^{\circ}(L)\big)^{\circ}\Big)-b\Big(\big(L-p^{\circ}(L)\big)^{\circ}\Big)\bigg)^{\circ}+a(L).$$

\end{proposition}

\vskip 3mm
\subsection{Proof of Theorem \ref{theo:nichtinjective}.}
The proof is a consequence of the  lemmas in this subsection.  The last one gives  the
result.
\vskip 2mm
\begin{lemma}\label{A-p-var}(\cite{MeyerSchuettWerner2012}, Lemma 6)
Let $p \in {\mathfrak P}_n$  and let $g$ be the centroid. For $0 < \varepsilon  < 1$, define 
$A_{p,\varepsilon},B_{p,\varepsilon}:\mathcal K_{n}\rightarrow\mathcal K_{n}$ by
$$
A_{p,\varepsilon}(K)
=\left\{x\in K\left|  \   \left\langle x,p\Big(\big(K-g(K)\big)^{\circ}\Big)\right\rangle 
\geq     \sup_{y\in K} \left\langle y,p\Big(\big(K-g(K)\big)^{\circ}\Big) \right\rangle
-\varepsilon \right.\right\}.
$$
and
$$
B_{p,\varepsilon}(K)
=\left\{x\in K\left| \ \left\langle x,p\Big(\big(K-g(K)\big)^{\circ}\Big)\right\rangle 
\leq     \inf_{y\in K}\left\langle y,p\Big(\big(K-g(K)\big)^{\circ}\Big)\right\rangle
+\varepsilon \right.\right\}.
$$
Then  $A_{p,\varepsilon}$ and $B_{p,\varepsilon}$ are affine invariant set maps.
\end{lemma}
\vskip 2mm

\noindent
\begin{remark}
Since $0$ is the Santal\'o point of 
$\big(K-g(K)\big)^{\circ}$, 
$0\in{\mathfrak P}_n\Big(\big(K-g(K)\big)^{\circ}\Big)$. Therefore ${\mathfrak P}_n\Big(\big(K-g(K)\big)^{\circ}\Big)$ is a subspace
of $\mathbb R^{n}$.
\end{remark}
\vskip 2mm

\vskip 3mm
\begin{lemma} Define $p_{\varepsilon,\delta}:\mathcal K_{n}\rightarrow
\mathbb R^{n}$  by
$$
p_{\varepsilon,\delta}(K)
=g(A_{g,\varepsilon}(K)\cup B_{g,\delta}(K)).
$$
Then $p_{\varepsilon,\delta}$ is a proper affine invariant point.
\end{lemma}
\vskip 2mm
\noindent
{\bf Proof.} The sets
$A_{g,\varepsilon}(K)$ and $ B_{g,\delta}(K)$ have non-empty interior. Therefore,
$p_{\varepsilon,\delta}(K)$ is well defined and it is an interior point of $K$.
By  Lemma \ref{A-p-var}, for every bijective, affine mapping $T:\R^n\to \R^n$ 
$$
A_{g,\varepsilon}(T(K))=T(A_{p,\varepsilon}(K))
\hskip 10mm
\mbox{and}
\hskip 10mm
B_{g,\delta}(T(K))
=T(B_{g,\delta}(K)).
$$
Therefore,
$$
A_{g,\varepsilon}(T(K))\cup B_{g,\delta}(T(K))
=T(A_{p,\varepsilon}(K)\cup B_{g,\delta}(K)).
$$
It follows  also from Lemma \ref{A-p-var} that $p_{\varepsilon,\delta}$ is continuous.  
$\Box$
\vskip 2mm

\noindent
\begin{remark} If $A_{g,\varepsilon}(K)\cap B_{g,\delta}(K)=\emptyset$,
then
\begin{eqnarray}\label{schwerpunkt11}
&& g(A_{g,\varepsilon}(K)\cup B_{g,\delta}(K)) \\
&&=\frac{|A_{g,\varepsilon}(K)|}{|A_{g,\varepsilon}(K)|+|B_{g,\delta}(K)|}
\  g(A_{g,\varepsilon}(K))
+
\frac{|B_{g,\delta}(K|}{|A_{p,\varepsilon}(K)|+|B_{g,\delta}(K)|}
\  g(B_{g,\delta}(K)). \nonumber
\end{eqnarray}
\end{remark}
\vskip 3mm
\noindent
\begin{lemma}\label{Trapez1}
For numbers $0<a<b$  let $K(a,b)$ be the convex body in $\R^2$ defined by
$$
K(a,b)
=\operatorname{conv}
\left\{\left(-\frac{\frac{2}{3}b+\frac{1}{3}a}{a+b},\pm a\right),
\left(\frac{\frac{2}{3}a+\frac{1}{3}b}{a+b},
\pm b\right)
\right\}.
$$
Then
$$
g\big(K(a,b)\big)=\left(0,0\right)\hbox{
\hskip 3mm and\hskip 3mm }
g\big(K(a,b)^{\circ}\big)
=\left(\frac{-3ab(b^2-a^2)}{(2a^2+2b^2+5ab)(2a^2+2b^2+2ab)}, \  0\right).
$$
In particular, the first coordinate of $g\big(K(a,b)^{\circ}\big)$ is negative.
\end{lemma}
\vskip 2mm

\begin{remark} Note  that $K(a,b)$ is the  translate of $ \operatorname{conv}\left\{\left(0,\pm a\right),
\left(1,\pm b\right)\right\}$ by $\left(-\frac{\frac{2}{3}b+\frac{1}{3}a}{a+b},  \   0\right)$.
\end{remark}
\vskip 2mm
\noindent
{\bf Proof.}
By symmetry the second coordinates of $g\big(K(a,b)\big)$
and $g\big(K(a,b)^{\circ}\big)$ are $0$. That $
g\big(K(a,b)\big)=\left(0,0\right)$ follows from a simple computation.
We see that 
\begin{eqnarray*}
K(a,b)=\left\{(x,y)\in \R^2:   -\frac{\frac{2}{3}b+\frac{1}{3}a}{a+b} \le x\le \frac{\frac{2}{3}a+\frac{1}{3}b}{a+b}, \   |y| \le (b-a) x+ \frac{2(a^2+b^2+ab)}{3(a+b)}\right\}.
\end{eqnarray*}
Hence
\begin{eqnarray*}
K(a,b)^{\circ}=\operatorname{conv}
\left\{     \left(-\frac{3(a+b)}{a +2b},0\right),
\left(\frac{3(a+b)}{2a+b},0\right),\left(\frac{-3(b^2-a^2)}{2a^2+2b^2+2ab},\frac{\pm 3(b+a)}{2a^2+2b^2+2ab}\right)     \right\}.
\end{eqnarray*}
Therefore,
\begin{eqnarray*}
g\big(K(a,b)^{\circ}\big)&=&
\left(\frac{1}{3} \left(-\frac{a+b}{\frac{2}{3}b+\frac{1}{3}a}+\frac{a+b}{\frac{2}{3}a+\frac{1}{3}b}-3\frac{b^2-a^2}{2a^2+2b^2+2ab}\right), \  0\right)\\
&=& \left(\frac{-3ab(b^2-a^2)}{(2a^2+2b^2+5ab)(2a^2+2b^2+2ab)}, \  0\right).\\
\end{eqnarray*}
In particular,  the first coordinate of
$g \big( K(a,b)^{\circ} \big) $ is negative. 
 $\Box$
\vskip 3mm

\begin{lemma}\label{Noninj1}
(i) For all $\varepsilon,\delta>0$,  we have  that 
$
p_{\varepsilon,\delta}(B_{\infty}^{2})=0.
$
\vskip 2mm
\noindent
(ii) For all $\eta\in (0,1)$ there exist $\varepsilon>0$ and $\delta>0$
 such that 
 $$
p_{\varepsilon,\delta}\big((B_{\infty}^{2})_{\eta e_{1}}\big)=0.
$$
Here, $e_1=(1,0)$.
\end{lemma}
\vskip 2mm
\noindent
{\bf Proof.}
(i) Since $B_{\infty}^{2}$ is $0$-symmetric,  we have that 
$
p_{\varepsilon,\delta}(B_{\infty}^{2})=0.
$
\par
\noindent
(ii) 
For $\eta$ with $0\leq\eta<1$, let $B_{\eta}:=(B_{\infty}^{2})_{\eta e_1}$. Then
$$
B_{\eta}=
\operatorname{conv}
\left\{\left(-\frac{1}{1+\eta}, \pm\frac{1}{1+\eta}\right),
\left(\frac{1}{1-\eta}, \pm\frac{1}{1-\eta}\right)\right\}.
$$
We put  $a=\frac{1}{1+\eta}$ and $b=\frac{1}{1-\eta}$ and  get that
$g(B_{\eta})= \big(\frac{2(b-a)}{3},0\big)$. Hence
$$B_{\eta}-g(B_{\eta})= \operatorname{conv}
\left\{
\left(-\left(\frac{2}{3}b +\frac{1}{3} a\right), \ \pm a\right), \left( \frac{2}{3}a +\frac{1}{3} b\  \pm b \right) \right\}.$$
Thus,  with the notation of Lemma \ref{Trapez1},
 $$B_{\eta}-g( B_{\eta} )=T_{a,b} \big( K(a,b)\big),$$
 where $T_{a,b}:\R^2 \to \R^2$ is defined by
$$T_{a,b} (x,y)=\big((a+b)x, y\big).$$
It follows  with (\ref{adjoint}) that
$$\big(B_{\eta}-g( B_{\eta} )\big)^{\circ} =(T_{a,b}^*)^{-1} \big(K(a,b)^{\circ}\big).$$
Again by Lemma \ref{Trapez1}, 
$$g\Big(\big(B_{\eta}-g( B_{\eta}) \big)^{\circ} \Big)= 
 \left(\frac{-3ab(b-a)}{(2a^2+2b^2+5ab)(2a^2+2b^2+2ab)}, \  0\right).$$
We replace $a=\frac{1}{1+\eta}$ and $b=\frac{1}{1-\eta}$ and set
$G(\eta):=g\Big(\big(B_{\eta}-g( B_{\eta}) \big)^{\circ} \Big)$.
Then
$$G(\eta)=\big(\alpha(\eta),0\big)= \left(\frac{-3 \eta (1-\eta^2)^2}{(3+\eta^2)(9-\eta^2)}, \  0 \right),$$
with $\alpha(\eta)<0$ for every $\eta\in (0,1)$.
Now we compute
$$p_{\varepsilon,\delta}(B_{\eta})
=
g\big(A_{g,\varepsilon}(B_{\eta})
\cup B_{g,\delta}(B_{\eta})\big).
$$
Since $\alpha(\eta)<0$ for $0<\eta<1$
\begin{eqnarray*}
A_{g,\varepsilon}(B_{\eta})  &=&\{x\in B_{\eta}|\ 
\langle x,G(\eta)\rangle
\geq\sup_{z\in B_{\eta}}\langle z,G(\eta)\rangle-\varepsilon\}
\\
&=&\{x\in B_{\eta}| \ 
\langle x,\alpha(\eta)e_{1}\rangle
\geq\sup_{z\in B_{\eta}}\langle z,\alpha(\eta)e_{1}\rangle-\varepsilon\}\\
&=&\left\{x\in B_{\eta}\left|\ 
\langle x,-e_{1}\rangle
\geq\sup_{z\in B_{\eta}}\langle z,-e_{1})\rangle-\frac{\varepsilon}{|\alpha(\eta)|}\right.\right\}   \\
&=&\left\{(x_{1},x_{2})\in B_{\eta}\left|\ 
-x_{1}
\geq\frac{1}{1+\eta}-\frac{\varepsilon}{|\alpha(\eta)|}\right.\right\}    \\
&=&\left\{(x_{1},x_{2})\in B_{\eta}\left|\ 
x_{1}
\leq-\frac{1}{1+\eta}+\frac{\varepsilon}{|\alpha(\eta)|}\right.\right\}.
\end{eqnarray*}
Similarly, 
\begin{eqnarray*}
B_{g,\delta}(B_{\eta})
&=&\left\{x\in B_{\eta}\left|\ 
x_{1}
\geq\frac{1}{1-\eta}-\frac{\delta}{|\alpha(\eta)|}\right.\right\}.
\end{eqnarray*}
Please note that
\begin{equation} \label{epsilon und delta}
A_{g,\varepsilon}(B_{\eta})
\cap
B_{g,\delta}(B_{\eta})
=\emptyset \hskip 3mm \text{if and only if} \hskip 3mm 
\varepsilon +\delta< \frac{2|\alpha(\eta)|}{1-\eta^2}.
\end{equation}
Suppose that we take $\eps$ and $\delta$ so that this condition holds.
By (\ref{schwerpunkt11}), 
$$p_{\varepsilon,\delta}(B_{\eta})=
\frac{|A_{g,\varepsilon}(B_{\eta})|   g\big( A_{g,\varepsilon}(B_{\eta})\big)+ 
        |B_{g,\delta}        (B_{\eta})|   g\big(B_{g,\delta}          (B_{\eta})\big)  }
       {  |A_{g,\varepsilon}  (B_{\eta} )| + |B_{g,\delta}  (B_{\eta} )|  }.$$
Thus $p_{\eps,\delta}(B_{\eta})=0$ if and only if 
\begin{eqnarray*}
0&=& |A_{g,\varepsilon}(B_{\eta})|   g\big( A_{g,\varepsilon}(B_{\eta})\big)+ 
        |B_{g,\delta}        (B_{\eta})|   g\big(B_{g,\delta}          (B_{\eta})\big).        
 \end{eqnarray*} 
Since $$B_{\eta}=\left\{(x,y)\in \R^2 \big| -\frac{1}{1+\eta} \le x\le \frac{1}{1-\eta},\ |y| \le 1+\eta x\right\},$$ 
and
$$
|A_{g,\varepsilon}(B_{\eta})|  = \frac{\varepsilon}{|\alpha(\eta)|} \left( \frac{2}{1+\eta} +\frac{\varepsilon \eta}{|\alpha(\eta)|}\right), \ \  |B_{g,\delta}(B_{\eta})|  = \frac{\delta}{|\alpha(\eta)|} \left( \frac{2}{1-\eta} -\frac{\delta \eta}{|\alpha(\eta)|}\right),
$$ 
$p_{\eps,\delta}(B_{\eta})=0$ if and only if      
 \begin{eqnarray*}  
    0    &=& \left(\frac{\varepsilon}{|\alpha(\eta)|} \left( \frac{2}{1+\eta} +\frac{\varepsilon \eta}{|\alpha(\eta)|}\right)\right) \  \int_{\left\{      (x,y)\in B_{\eta}:
-\frac{1}{1+\eta} \leq x \leq -\frac{1}{1+\eta} +\frac{\varepsilon} {|\alpha(\eta)| }\right\} } xdxdy \\
&+& \left(  \frac{\delta}{|\alpha(\eta)|} \left( \frac{2}{1-\eta} -\frac{\delta \eta}{|\alpha(\eta)|}\right) \right) \  
\int_{\left\{(x,y)\in B_{\eta}:
 \frac{1}{1-\eta} -\frac{\delta} {|\alpha(\eta)| }    \leq x\leq \frac{1}{1-\eta} \right\} } xdxdy, 
\end{eqnarray*}
which is equivalent to 
 \begin{eqnarray*}  
    0    &=& \varepsilon^2  \left( \frac{2}{1+\eta} +\frac{\varepsilon \eta}{|\alpha(\eta)|}\right)\   
 \left(- \frac{1}{(1+\eta)^2} + 
 \frac{\varepsilon}{2 |\alpha(\eta)|} 
 \left( \frac{1-\eta}{1+\eta} + \frac{2 \varepsilon \eta}{3|\alpha(\eta)|} \right)
  \right) \\
&+&  \delta^2  \left( \frac{2}{1-\eta} -\frac{\delta \eta}{|\alpha(\eta)|}\right)\   
 \left( \frac{1}{(1-\eta)^2} - 
 \frac{\delta}{2 |\alpha(\eta)|} 
 \left( \frac{1+\eta}{1-\eta} - \frac{2 \delta \eta}{3|\alpha(\eta)|} \right)
  \right).
\end{eqnarray*}
For $\varepsilon \leq  \varepsilon_0 < \frac{|\alpha(\eta)|}{ 1 - \eta^2}$ fixed, put
\begin{eqnarray*} 
f_\varepsilon(\delta)  &=& \varepsilon^2  \left( \frac{2}{1+\eta} +\frac{\varepsilon \eta}{|\alpha(\eta)|}\right)\   
 \left(- \frac{1}{(1+\eta)^2} + 
 \frac{\varepsilon}{2 |\alpha(\eta)|} 
 \left( \frac{1-\eta}{1+\eta} + \frac{2 \varepsilon \eta}{3|\alpha(\eta)|} \right)
  \right) \\
&&+  \delta^2  \left( \frac{2}{1-\eta} -\frac{\delta \eta}{|\alpha(\eta)|}\right)\   
 \left( \frac{1}{(1-\eta)^2} - 
 \frac{\delta}{2 |\alpha(\eta)|} 
 \left( \frac{1+\eta}{1-\eta} - \frac{2 \delta \eta}{3|\alpha(\eta)|} \right)
  \right).
\end{eqnarray*}
Then, for sufficiently small $\epsilon$
$$
f_\varepsilon(\varepsilon) >0  \hskip 5mm \text{and} \hskip 5mm f_\varepsilon(\varepsilon^2) < 0.
$$
Thus, for all $\varepsilon \leq  \varepsilon_0$, by the Intermediate Value Theorem, there exists $\varepsilon^2  < \delta < \varepsilon$  with $f_\varepsilon(\delta) =0$, and hence 
$p_{\eps,\delta}(B_{\eta})=0.$
\hskip 2mm
$\Box$
\vskip 3mm\noindent
{\bf Proof of Theorem \ref{theo:nichtinjective}. } 
Lemma \ref{Noninj1} provides an example in dimension $2$. This example is easily generalized
to dimension $n$, with
$B_{\infty}^n$ instead of $B_{\infty}^2$. 
$\Box$

\vskip 3mm

\section{Examples of affine invariant points and sets.} \label{example}

In this Section, we list some of the  classical affine invariant points and sets, with proof if necessary. 
We will also introduce several  new affine invariant points and sets.
First, we state a  lemma (the proof of which we leave to the reader)  that provides a general tool to study  those affine points which are given
as minima or maxima of functions. 
\par
There, ${\cal C}_0(X)$ denotes  the space of continuous functions on  a locally compact metric space $X$, vanishing at $\infty$, endowed with the uniform norm $\|\ \cdot\|_\infty$.
\vskip 2mm\noindent
\begin{lemma}\label{functions} 
Let $\left(f_n\right)_{n \in \mathbb{N}} \subset {\cal C}_0(X)$ be a sequence of positive functions. Moreover assume that   $\|f_n-f\|_{\infty}\to 0$, where $f \in {\cal C}_0(X)$  reaches its maximum at a unique point $x\in X$.
Then for any sequence $(x_n)_{n\in \mathbb{N}}$ in $X$ such that
$f_n$ reaches its maximum at $x_n$ for all $n$, one has $x_n\to x$ in $X$. 
\end{lemma} 

\vskip 3mm
{\bf The John regions of a convex body.}
\vskip 2mm
Let $K\in \K$.  Let ${\cal E}_0$ be the set of all ellipsoids in ${\mathbb R}^n$
centered at $0$. We define a function
$f_K:\inte(K)\to {\mathbb R}_+$ by
$$f_K(x)=\sup\{|E| :  E\in {\cal E}_0,\ x+ E\subset K\}.$$
It is easy to  see that this supremum  is  a maximum, that $f_K$ is continuous on 
$\inte(K)$  and that $f_K(x)\to 0$ when  $x\to \partial K$. Thus $f_K$ can be extended to the whole ${\mathbb R}^n$ as a continuous function with compact support, setting 
$$f_K(x)=0\hbox{ for } x\in {\mathbb R}^n\setminus \inte(K).$$
\vskip 3mm
We omit the proof of the following easy lemma.
\vskip 2mm
\begin{lemma}\label{John1}
The mapping $K\to f_K$ is continuous from ${\mathcal K}_n$ to ${\cal C}_0( {\mathbb R}^n)$. 
\end{lemma}
\vskip 3mm
Let us recall the celebrated  theorem by F. John (see e.g.  \cite{NicoleBuch}).
\vskip 2mm
\begin{theorem}{(\bf F. John)} 
Let $K\in \K$ and  suppose that $B_2^n\subset K$ (resp. $K\subset B_2^n$). The following are equivalent. 
\par
\noindent
(i)  $B_2^n$ is the ellipsoid of maximal volume contained in $K$ (resp. of minimal volume containing $K$).
\par
\noindent
(ii) There exist  $u_i\in  S^{n-1}\cap\partial K\cap \partial K^*$ and    $c_i>0$, $1\le i\le m\le n(n+1)\ $,  such that
$$ \sum_{i=1}^m c_i u_i=0  \  \hbox{  and} \  \  
x=\sum_{i=1}^m c_i\langle x,u_i\rangle u_i \  \  \   \hbox{   for every}  \   x\in \R^n. $$
\end{theorem}
\vskip 3mm
Thus,  there is a unique ellipsoid of maximal volume  $J(K)$ contained in $K$, called the John ellipsoid of $K$, and a  unique ellipsoid of minimal volume  $L(K)$ containing  $K$, called the L\"owner ellipsoid of $K$.
We call its centers respectively, $j(K)$,  the {\it John point} of $K$ and $l(K)$, the L\"owner point of $K$. With the previous notation, $j(K)$ is the unique point $x\in \inte(K)$ such that $\|f_K\|_{\infty}=f_K(x)$.
\vskip 3mm
Then the next proposition follows immediately  from Lemmas \ref{functions}  and \ref{John1}.
\vskip 2mm
\begin{proposition}\label{John2}
 $K\to j(K)$ is an affine invariant point.
\end{proposition}
\vskip 3mm
The following lemma allows to say more.
\vskip 2mm
\begin{lemma}\label{konvex}
With the preceding notations, $f_K^{\frac{1}{n}}$ is concave on $\inte(K)$ and  hence $f_K$ is log-concave on ${\mathbb R}^n$.
\end{lemma}
\noindent
{\bf Proof.}   For $i=1,2$, let $x_i\in \inte(K)$ and $E_i\in {\cal E}_0$ such that 
$x_i+E_i\subset K$. After an affine transform, we may suppose that $E_1=B_2^n$ and $E_2=\Lambda  (B_2^n)$, where $\Lambda:{\mathbb R}^n\to {\mathbb R}^n$ is a diagonal matrix with positive entries $\lambda_1,\dots, \lambda_n$ on the diagonal. 
\par
For $t\in [0,1]$, define $E_t=\big((1-t) Id+t \Lambda\big)(B_2^n)$, where $Id$ is the 
identity matrix on ${\mathbb R}^n$. 
Then $E_t\in {\cal E}_0$ and 
$$(1-t)x_1+tx_2+E_t\subset (1-t)(x_1+E_1)+t(x_2+E_2)\subset K.$$
By the Brunn-Minkowski inequality (e.g., \cite{GardnerBook}, \cite{SchneiderBuch}), we get
\begin{eqnarray*}  
f_K^\frac{1}{n}\big(tx_1+(1-t)x_2\big)&\geq &|E_t| ^\frac{1}{n} =|(1-t) E_1+t E_2 |^\frac{1}{n}
\geq (1-t) |E_1|^{\frac{1}{n}}+t |E_2|^{\frac{1}{n}}.  \ \   \Box
 \end{eqnarray*}
 \vskip 3mm
 \noindent
\begin{definition} Let $c\in [0,1)$. We define the {\it John region of $K$ of index $c$} by
$$J_c(K)=\left\{x\in {{\mathbb R}}^n: f_K(x)\ge c \|f_K\|_{\infty}\right\}.
$$
\end{definition}

By  Lemma \ref{konvex}, $J_c(K)$ is convex.
With Lemma \ref{John1}, we then get the next proposition, which provides a new affine invariant set mapping.
\vskip 2mm
\noindent
\begin{proposition} For $ c\in (0,1)$, 
the mapping $K\to J_c(K)$ is a proper  affine  invariant  set mapping from
${\mathcal K}_n$ to ${\mathcal K}_n$.
\end{proposition}
\vskip 3mm
{\bf The L\"owner regions of a convex body.}
\vskip 1mm\noindent
Let $K\in \K$.
Define $\lambda_K: {\mathbb R}^n\to {\mathbb R}$ by 
$$\lambda_K(x)=\big( \inf\{ \vol_n(E):  E\in {\cal E}_0; K\subset x+E\}\big)^{-1}.$$
It is clear that $\lambda_K>0$ and $\lambda_K\in{\cal C}_0(\R^n)$.   The following lemma is easy.
\vskip 2mm\noindent
\begin{lemma}\label{loe}
 The mapping $K\to \lambda_K$ is continuous from ${\mathcal K}_n$ to ${\cal C}_0({\mathbb R}^n)$.
\end{lemma}
\vskip 2mm
Since the center $l(K)$ of the L\"owner ellipsoid $L(K)$ of $K$ is the unique point $x\in \inte(K)$ such that $\|\lambda_K\|_{\infty}=\lambda_K(x)$,  Lemmas \ref{functions} and \ref{loe} give the next proposition.
\vskip 2mm\noindent
\begin{proposition}
 $K\to l(K)$ is an affine invariant point.
\end{proposition}
\vskip 2mm
\noindent
\begin{definition}
If $c\in [0,1)$, we define the {\it L\"owner region of $K$ of index $c$} by
$$L_c(K)=\conv\big[\{x\in {{\mathbb R}}^n: \lambda_K(x)\ge c\|\lambda_K\|_{\infty}\}\big].$$
\end{definition}
\vskip 2mm

\begin{proposition}
$K\to L_c(K)$ is an affine invariant  set mapping from $\mathcal K_n$ to ${\mathcal K}_n$.
\end{proposition}
\vskip 2mm
In conclusion, Theorem \ref{Satz5} summarizes all these facts.
\vskip 2mm
\begin{theorem} \label{Satz5} Let  $K\in \K$.
If $j(K)$ and $l(K)$ denote respectively the centers of the John and of the L\"owner ellipsoids of $K$, then $l$ and $j$ are in $ \P$ and $l=j^{\circ}$.
\end{theorem}
\vskip 2mm
\noindent
{\bf Proof.} We only need  to prove that $l(K^{j(K)})=j(K)$, or, equivalently,  that if the John ellipsoid 
$E_K$
of $K$ is centered at $0$, then the  L\"owner ellipsoid of $K^{\circ}$ is 
$(E_K)^{\circ}$. This follows from John's theorem. \hskip 2mm $\Box$
\vskip 2mm
\noindent
\begin{remark} Note that we need the full strength of John's theorem only to prove that  $l=j^{\circ}$. The fact that $j$ and $l$ are uniquely defined follows from more elementary reasonnings.
\end{remark}

\vskip 3mm
{\bf The Santal\'o point and the center of gravity.}
 \vskip 2mm \noindent
The following result is well known (see e.g., \cite{Santalo}).
There, $S_K:\inte(K)\to \R^n$ is the function defined by
$S_K(x)=|K^x|$.
\vskip 2mm
\begin{theorem} \cite{Santalo}
Let $K\in \K$. Then the function $S_K$
is strictly log convex. Moreover   $S_K(x)\to +\infty$ when $x\to \partial K$ and $S_K$ reaches its minimum at a unique point $s(K)\in \inte(K)$. This point  is characterized
by the fact that $s(K)$ is the centroid of $K^{s(K)}$ (or that $0$ is the centroid of
$\big(K-s(K)\big)^{\circ}$).
\end{theorem}
\vskip 2mm\noindent
\begin{proposition}\label{santalo}
The mapping $K\to s(K)$ is a proper affine invariant point and 
$g=s^{\circ}$.
\end{proposition}
\vskip 2mm
\noindent
{\bf Proof.} The uniqueness of $s(K)$ shows that $g$ is injective. 
Hence, by Theorem \ref{theo:dual},  $g$ has  a dual point $g^{\circ}\in \P$.  Thus $g^{\circ}(K^{g(K)}) =g(K)$, 
 and by the preceding characterization, $s(K^{g(K)})=g(K)$.  Since $g$ is surjective by Theorem \ref{theo:surjective}, it follows that $g^{\circ}=s \in \P$. \hskip 2mm $\Box$
\vskip 2mm
\noindent
\begin{remark}  (i) The  fact that the mapping $\psi_K:\inte(K)\to \R^n$, $\psi_K(x)= g(K^x)\hbox{ for }x\in \inte(K)$, 
is  bijective can also be proved in an other way: 
The function $\Theta_K: \inte(K)\to (0,+\infty) $, 
$\Theta_K(x) =\log S_K(x)$,  is strictly convex and $\Theta(x)\to+\infty$ when $x\to \partial K$. It follows that $\nabla \Theta_K:\inte(K)\to \R^n$ is bijective. Moreover,  it is easily checked that
$$(\nabla\Theta_K)(x)=g\big((K-x)^{\circ}\big)\hbox{ for all } x\in \inte(K).$$
\par
\noindent
(ii) The Santal\'o regions of $K$,  defined in \cite{MeyerWerner1998} , for $c>0$ by 
$$
S_c(K) = \left\{z \in \inte(K): \left|K^z \right| \leq (1+c) \   \left|K^{s(K)}\right|\right\}
$$ 
are affine invariant set mappings.
\end{remark}
 \vskip 3mm
{\bf The center of the maximal  volume centrally symmetric body inside $K$.} 
  \vskip 1mm \noindent
The first part of the following proposition follows from the Brunn-Minkowski inequality, together
with its equality case, the second part   from Lemma  \ref{functions}.
\vskip 2mm
\begin{proposition}\label{maxvol}
Let $K\in {\mathcal K}_n$. Then the function $\theta_K(x) =\vol_n\big(K\cap(2x-K)\big)^\frac{1}{n}$  
is concave on its support and reaches its maximum at a unique point $m(K)$. 
Moreover, the mapping $m$ is  a proper  affine invariant point.
\end{proposition}
\vskip 1mm \noindent
\begin{proposition}
Let $0<c<1$.  For $K\in \K$, define
$$ M_c(K) = \{x\in  {\mathbb R}^n: |K\cap(2x-K)|\ge  c|K\cap(2m(K)-K)|\}.$$
Then $K\to M_c(K)$ is an affine invariant set mapping.
\end{proposition}
\vskip 1mm
\noindent
{\bf Proof.}  As $\theta_K^n$  is concave on its support, $M_c(K)$ is convex. Affine invariance and continuity of the map
$K\to M_c(K)$ are easy.   $\Box$
 \vskip 3mm 
{\bf The center of the maximal  volume zonoid body inside $K$.} 
  \vskip 1mm \noindent
Let ${\cal Z} $ be the (closed) set of zonoids in $\K$, and ${\cal Z} _0$ be the set of all zonoids  that are centered at $0$ (see\cite{GardnerBook} or \cite{SchneiderBuch}).  For $K\in \K$, let
\begin{eqnarray*}
g_K(x) = \left\{
\begin{array}{cc}  \max\{|Z|: Z\in{\cal Z} _0,\  x+Z\subset K\}   &  \text{if }  x\in \inte(K)\\
 0   & \text{if }  x\not\in \inte(K).
 \end{array}
 \right.
\end{eqnarray*}
It is clear that $K\to g_K$ is continuous  from $\K$ to $C_0(\R^n)$.
Since convex combinations (for the Minkowski addition) of zonoids are zonoids, it follows 
 as in Proposition \ref{maxvol},   that  $g_K$ reaches its maximum at a unique point $z(K)$ and that  $g_K^{1/n}$ is concave on its support. Thus
$$K_z(c)=\{x\in \R^n: g_K(x)\ge c\|g_K\|_\infty\}, \hskip 3mm 0\le c\le 1, $$
is convex.  And, again with Lemma \ref{functions}, we get the following proposition.
\vskip 2mm
\begin{proposition}\label{zcenter} The mapping $z: \K \to \R^n$ is a proper  affine invariant point and,  for every $0< c\le 1$, $K\to K_z(c)$ is a proper affine invariant set mapping. 
\end{proposition}
\vskip 2mm 

\noindent
\begin{remark} 
(i)  For $K\in \K$  and  $x\in \inte(K)$, define
$$\phi_K(x)=\max\{\lambda>0: \ x-\lambda\big(K-g(K)\big)\subset K - g(K)\},$$
where $g(K)$ is the centroid of $K$. Then $\phi_K$ is positive and concave
on $\inte(K)$ and $\phi_K(x)\to 0$
when  $x\to \partial K$. So one can extend $\phi_K$ to a
continuous
function on $\R^n$ by setting $\phi_K(x)=0$ when $x\not\in \inte(K)$.  It
is well known that
$\max_{x} \phi_K(x)\ge \phi_K(g(K)) \ge \frac{1}{n}\ . $ It follows easily
that for any $0<\delta<1$,
$$
K\to \{ \phi_K \ge (1-\delta)\max \phi_K\}
$$
is an affine invariant set mapping.
\par
But generally, $\phi_K$ does not reach its maximum at a unique point.  To
see that,
take $K=\Delta_2\times [-1,1]^{n-2} \subset \R^n$, where $\Delta_2$ is a
regular simplex centered at  $0$ in $\R^2$.
Then
$$
\{x \in \mathbb{R}^n: \phi_K(x)=\max \phi_K\} =\{(0,0,x_3,\dots, x_{n}):  |x_i|\le\frac{1}{2},  \  3\le i \le n \}.$$
\par 
(ii) Instead of the centers of centrally symmetric convex bodies contained in  $K$, we may study  the centers of those containing $K$.
For $K\in {\mathcal K}_n$, we define a positive and convex function
$\rho_K: {\mathbb R}^n\to {\mathbb R}_+$ by 
$$\rho_K(x)= |{\rm conv}[K,2x-K]|.$$
 If $L\in \K$ is  centered at $x\in \R^n$ and satisfies $K\subseteq L$, then $2x-K\subseteq L$ and therefore $\conv(K,2x-K)\subseteq L$. It follows that
$$\min\{|L|: L\in \K, K\subseteq L,  L\hbox{ is centrally symmetric}\}=\min_x \rho_K(x).$$
It is clear that $\rho_K\to +\infty$,  when $|x|\to \infty$. For
 $c>1$, let
$$N_c(K)=\{x\in {\mathbb R}^n: \rho_K(x)\le c \min_{z\in {\mathbb R}^n}
\rho_K(z)\}.$$ 
Then $K\to N_c(K)$ is an affine invariant  set mapping. 
\par
However,  one cannot define an affine invariant point in that
way, because it may happen that $\rho_K$  does not reach its
minimum at a unique point. 
For instance, if $K$ is a simplex in $\R^n$ and 
$n$ is even, then
$\{\rho_K =\min \rho_K\}$ has non-empty interior.
\end{remark}

\vskip 3mm 
{\bf  The illumination body.} 
 \vskip 1mm 
 \noindent
Let $K\in{\mathcal K}_n$, $\delta\ge 0$, $x\in \R^n$ and $F_K(x)= |\conv(x, K)|$. 
The {\it illumination body $K^{\delta}$}  \cite{Werner1994} of $K$ is defined by
$$
K^{\delta} = \left\{x\in {\mathbb R}^n: F_K(x)\le(1+\delta)|K| \right\}.
$$
Then  $K^{\delta}\in {\mathcal K}_n$. Clearly $K^0=K\subseteq K^{\delta}$ and 
$K\to K^{\delta}$ is affine invariant.

If $\sigma_K$ is the surface measure of $K$ on $S^{n-1}$ and $h_K :\R^{n} 
\to \R$,  $h_{K}(\xi)=\sup_{x\in K} \ \langle \xi,x \rangle$, 
is the support function of $K$, then 
$$F_K(x)= \frac{1}{n}\int_{S^{n-1}} \max\big(\langle x,u\rangle, h_K(u)\big)\ d\sigma_K(u).$$
Thus $F_K$ is continuous, convex and clearly $F_K(x)\to +\infty$, when $|x| \to +\infty$. 
\vskip 2mm
\noindent
\begin{proposition}\label{ill}
Let $I^{\delta}:\K\to \K$ be defined by $I^{\delta}(K)=K^{\delta}$. Then $I^{\delta}\in \A$. 
\end{proposition}
\noindent
{\bf Proof.}  We prove the continuity. Fix $\delta>0$ and $K\in \K$. We may suppose that $g(K)=0$,   and that for some $0 < r \leq R < \infty$, $rB_2^n\subseteq K\subseteq RB_2^n$. For $0 < \eta <1$, there exists $\eps>0$, such that  for all $L\in \K$ with $d_H(K,L)\le \eps$, one has,    $(1-\eta)K\subseteq L\subset (1+\eta)K$. Thus, for $x\in \R^n$, 
\begin{eqnarray*}         
    (1-\eta) \ \conv \big[x,K \big] &\subset& \conv\big[(1-\eta)x, (1-\eta)K \big] \subseteq \conv\big[x, (1-\eta)K\big] 
         \subset \conv\big[x,L\big]\\
         &\subset&
          \conv\big[x, (1+\eta)K\big] \subset \conv\big[(1+\eta)x, (1+\eta)K\big] \\
          &=&(1+\eta)\conv\big[x,K\big].
          \end{eqnarray*}
It follows that for $t>0$,    
$$\left\{x\in {\mathbb R}^n: F_K\le \frac{t}{(1+\eta)^n}\right\}\subseteq \{x\in {\mathbb R}^n: F_L\le t\} \subseteq  \left\{x\in {\mathbb R}^n: F_K\le \frac{t}{(1-\eta)^n}\right\}.$$        
Moreover,  $|\,|K|-|L|\,|  \le \rho$,   where $\rho = 2 n \eta \max\{|K|,|L|\}$. Consequently, 
 $$\left\{x\in {\mathbb R}^n: F_K\le \frac{(1+\delta)(1-\rho)}{(1+\eta)^n} |K|\right\}\subseteq L^{\delta}\subseteq \left\{x\in {\mathbb R}^n: F_K\le \frac{(1+\delta)(1+\rho)}{(1-\eta)^n}|K|\right\}.$$
 This allows to conclude,  because $\delta\to \{F_K\le (1+\delta)|K|\}$ is continuous. \hskip 2mm $\Box$
\vskip 3mm
{\bf The convex floating body.}
\vskip 1mm
\noindent
For a convex body $K$ and for $0\le \delta <\left(\frac{n}{n+1}\right)^n$, the convex floating body $K_{\delta}$ of $K$ was defined in \cite{SchuettWerner1990}
as the intersection  of all halfspaces $H^+$ whose
defining hyperplanes $H$ cut off a set of volume at most $\delta |K|$
from $K$, 
\begin{equation*}
K_{\delta}=\bigcap_{|H^-\cap K| \leq \delta  |K|} {H^+}. 
\end{equation*}
Then  the map $F_{\delta}:\K\to \K$, defined by $F_{\delta}(K)=K_{\delta}$ is an affine invariant set mapping.
This, and affine invariant points defined via the convex floating body are treated  in \cite{MeyerSchuettWerner2012}.
\vskip 3mm 
{\bf Extension to subsets of $\R^n$ of  affine invariant points on ${\mathcal K}_{k}$, $1\le k\le n-1$.} 
\vskip 2mm \noindent
Let ${\K}_{,k}$  be the set of closed  convex subsets $L$ of ${\mathbb R}^n$, 
whose affine span $E_L$  is $k$-dimensional. If $L\in {\K}_{,k}$, there exists a (non unique) one-to-one affine map $U: {\mathbb R}^n={\mathbb R}^{k}\times {\mathbb R}^{n-k}\to {\mathbb R}^n$
such that $U({\mathbb R}^{k})=E_L$. For $p\in {\mathfrak P}_{k}$, we then define 
\begin{equation}\label{def:k-dim p}
p(L)= U\Big(p\big(U^{-1}(L)\big)\Big).
\end{equation}
It is easy to show that this definition does not depend of the choice of $U$. Moreover this extended $p$ satisfies $p(L)\in E_L$ for every $L\in {\K}_{,k}$
and it is {\it affine invariant}:  For any one-to-one  affine mapping $W:{\mathbb R}^n\to {\mathbb R}^n$, one has
 $$p\big((W(L)\big)=W\big(p(L)\big).$$
\vskip 2mm
\begin{definition}\label{extaip}
We then call  $p : {\K}_{,k}\to \R^n $ defined by  (\ref{def:k-dim p}),  an extended affine invariant point.
\end{definition}
\par
\noindent
\begin{example} 
For $k=1$, the unique invariant point is the midpoint of a segment. It is canonically extended to be the midpoint of any segment of $\R^n$.  Another natural  example is the centroid of a body in $\R^k$, which extends  to the centroid of $L$ in $E_L$, when $L\in {\K}_{,k}$.
\end{example}
\vskip 3mm
The following proposition summarizes all these facts.
\vskip 2mm
\begin{proposition}\label{extension}
Let $p:{\mathcal K}_{n,k}\to\R^n$  be the affine invariant  extension of $p: {\mathcal K}_{k}  \to \R^k $. Then
$p(L)\in E_L$ for any $L \in {\mathcal K}_{n,k}$. $p$ is continuous, when ${\mathcal K}_{n,k}$ is endowed with the Hausdorff metric. Moreover,  this extension is proper (in the sense that $p(L)$ is in the relative interior of $L$ in $E_L$),  whenever $p:{\mathcal K}_{k}\to\R^k$ is  proper.
\end{proposition}
\vskip 3mm 
Let ${\cal M}_+({\mathbb R}^n)$ be the cone of non-negative Radon measures on ${\mathbb R}^n$.
 For a Borel function $f: \mathbb {R}^n \rightarrow \mathbb {R}^n$ and $\mu \in {\cal M}_+({\mathbb R}^n)$, let  $f(\mu)$ be the image measure of $\mu$ by $f$, i.e.
$$(f(\mu))(B)=\mu(f^{-1}(B)), \hbox{ for any Borel subset $ B$ of ${\mathbb R}^n$}.$$
For a  map $M: {\mathcal K}_n\to {\cal M}_+({\mathbb R}^n)$,  denote $\mu_K = M(K)$.
\vskip 3mm
\begin{definition}
We say that a mapping $M:\K \to {\cal M}_+({\mathbb R}^n)$ is a an affine invariant measure map if it is continuous for the Hausdorff topology on ${\cal K}_n$ and the weak*-topology on ${\cal M}_+({\mathbb R}^n)$ and if  
$\mu_{UK}=U(\mu_K)$ for any affine one-to-one mapping $U:{\mathbb R}^n\to {\mathbb R}^n$. Thus
 $$\int h(y)d\mu_{UK}(y)= \int h(Ux)d\mu_K(x), \  \hbox{ for any non-negative Borel function $h$ on $\R^n$}.$$
\end{definition}

\vskip 2mm\noindent
\begin{example} (i)  Let $A\in {\mathfrak A}_n$.
Let $\mu_K= \frac{{\bf 1}_{A(K)}}{|A(K)|}  dx$.  Then
 for every non negative Borel function $h$, 
$$\int h(x) d\mu_K(x)=  \frac{1}{|A(K)|} \int_{A(K)} h(x) dx.$$
If $U:{\mathbb R}^n\to {\mathbb R}^n$ is  an affine one-to-one mapping, then
\begin{eqnarray*}
\int h(y)d\mu_{UK}(y)&=& \frac{1}{|A(UK)|}\int_{A(UK)} h(y) dy
= \frac{1}{|(U(A(K))|}\int_{(U(A(K))} h(y) dy\\
&=&\frac{1}{|\det(U)| \ |A(K)|} 
\int_{A(K)} |\det(U)|h(Ux) dx= \int h(Ux)d\mu_{K}(x).
\end{eqnarray*}
\par
(ii) Let $K\in \K$  and let $p\in {\mathfrak P}_n$  be proper. For $x\in {\mathbb R}^n$, let
$$
\|x-p(K)\|_{K-p(K)}=\inf\{\lambda\ge 0: x-p(K)\in \lambda\big(K-p(K)\big)\}.
$$
If $\phi:{\mathbb R}_+\to {\mathbb R}_+$ is a Borel function, let 
\begin{equation}\label{fi}
\mu_K =\frac{1}{|K|}\phi\big(\|x-p(K)\|_{K-p(K)}\big) dx.
\end{equation}
Then, as above,  it is easy to see that $K\to \mu_K$ is an affine invariant measure map.
\par
(iii) Let  $a$ and $\delta$ be strictly positive real numbers.
We take $\phi(t)=\frac{1}{\delta} {\bf 1}_{[a,a+\delta)}$ in (\ref{fi}) of Example (ii). Let $L=K-g(K)$. Then for any Borel function $h$, when $\delta\to 0$, 
\begin{eqnarray*}
\int h(x) d\mu_K(x) &=& \frac{1}{\delta|L|}\int_{a\le \|x-g(K)\|_L\le a+\delta} h(x) dx \\
&=&\frac{n|B^n_2|}{\delta|L|}\int_{S^{n-1}}  \left( 
\int_{\frac{a}{ \| \theta\|_L } \le r\le \frac{a+\delta} {\|\theta\|_L } } 
 h(r\theta) r^{n-1} dr\right)  d\sigma(\theta)\\
&\sim& \frac{ v_n}{|L|} \int_{S^{n-1}}
h\left( \frac{a\theta}{||\theta||_L}\right )   \frac{1} {\|\theta\|_L^n} d\sigma(\theta).
\end{eqnarray*}
Here,  $\sigma$ is the normalized measure on $S^{n-1}$. Therefore, 
$K\to \mu_K$ is  an affine invariant measure map, with $\mu_K$  given by
$$\int h(x) d\mu_K(x) =\frac{1}{|L|}\int_{S^{n-1}} h\left( \frac{a\theta}{||\theta||_L}\right )   \frac{1} {\|\theta\|_L^n} d\sigma(\theta)
.$$
\end{example}

\vskip 3mm
With the previous definition, the following result is easy to prove.
\vskip 2mm
\begin{proposition} Let 
 $p\in {\mathfrak P}_n$, let $q: {\K}_{,n-1}\to {\mathbb R}^n$
be an extended affine invariant
point  and  let $\mu$
 be an affine invariant measure such that  $\mu_K=\mu(K)$ is supported by
$K$ for all $K\in\K$. Then the map $d:{\mathcal K}_{n}\to {\mathbb R}^n$ defined by 
$$K\to d(K) = p(K)+\int q\big( \{x\in K-p(K): \langle x^*,x\rangle =1\}\big)d\mu_{\big(K-p(K)\big)^\circ}(x^*), $$
is an affine invariant point.
\end{proposition}
\vskip 2mm
\noindent
\begin{example} For $0<\alpha<1$ and $w\in \P$, 
we define $d_1$ and $d_2\in  {\mathfrak P}_n$ by
$$ 
d_1(K)=p(K)+ \frac{1}{\left|\big(K-p(K)\big)^\circ\right|} \int_{S^{n-1}} \frac{q\left(\left\{x\in K-p(K): \langle \theta, x\rangle =\alpha\|\theta\|_{\left(K-p(K)\right)^\circ}\right\}\right) }{\|\theta\|_{\left(K-p(K)\right)^\circ}^n} d\sigma(\theta).
$$ 
$$ 
d_2(K)=p(K)+ \frac{1}{\left|\big(K-p(K)\big)^\circ\right|} \int_{S^{n-1}} \frac{w\left(\left\{x\in K-p(K): \langle \theta, x\rangle \ge \alpha\|\theta\|_{\left(K-p(K)\right)^\circ}\right\}\right) }{\|\theta\|_{\left(K-p(K)\right)^\circ}^n} d\sigma(\theta).
$$
\end{example}

\vskip 3mm
{\bf Two affine invariant points related to the projection body.}
\vskip 1 mm\noindent
Let $K\in \mathcal{K}_n$. We recall that the support function $h_{\Pi K}:\R^n \to \R_+$ of the projection body $\Pi K$ of $K$  (see e.g., \cite{GardnerBook, 
SchneiderBuch}) is given  by \begin{eqnarray*}
h_{\Pi K}(x)= \left\{
\begin{array}{cc} | x | P_{u_x}K|  &  \text{if }  x\not=0\\
 0 & \text{ if }  x=0,
 \end{array}
 \right.
\end{eqnarray*}
where $u_x=\frac{x}{|x|} \in S^{n-1}$ and for $u\in S^{n-1}$, $P_u:\R^n\to \R^n$  is the orthogonal projection onto $\{u\}^{\perp}$.
It is well known (see \cite{GardnerBook, SchneiderBuch}) that for any one-to-one linear map $T:\R^n\to \R^n$, one has
$$h_{\Pi(TK)}(Tx)= |\det(T)|\ h_{\Pi K}(x) \hbox{ for all } x\in \R^n.$$
\par
We define two affine invariant points related to the projection body. In both definitions we use the centroid $g$ which
could be replaced by any other $p\in {\mathfrak P}_n$.
\vskip 2mm
(i) For $K\in\K$, let
$\pi(K)= \frac{1}{|K|^2} \int_{K-g(K)} x\ h_{\Pi K}(x) dx +g(K).$
Then $\pi:\K\to \R^n$ is  an affine invariant point. Observe that 
$$\pi(K)-g(K)=\frac{1}{|K|^2} \int_{K-g(K)} x\ h_{\Pi K}(x) dx=\frac{1}{(n+2) |K|^2} \int_{S^{n-1}} \frac{\theta}{\|\theta\|_{K-g(K)}^{n+2}} h_{\Pi K}(\theta) d\sigma(\theta).$$
\vskip 2mm
(ii)  For $K\in\K$, let $b_t(K)=\frac{1}{|K|} \int_{K-g(K)} g(K+t[-x,x]) dx.$

\vskip 4mm

Mathieu Meyer \\
  {\small        Universit\'{e} de Paris Est - Marne-la-Vall\'{e}e}\\
    {\small      Equipe d'Analyse et de Math\'ematiques Appliqu\'{e}es}\\
   {\small       Cit\'e  Descartes - 5, bd Descartes }\\
    {\small      Champs-sur-Marne
          77454 Marne-la-Vall\'{e}e,  France} \\
          {\small \tt mathieu.meyer@univ-mlv.fr} \\

   Carsten Sch\"utt \\
     {\small       Christian Albrechts Universit\"at }\\
        {\small    Mathematisches Seminar }\\
        {\small    24098 Kiel, Germany} \\
         {\small \tt schuett@math.uni-kiel.de}   \\

     \and Elisabeth Werner\\
{\small Department of Mathematics \ \ \ \ \ \ \ \ \ \ \ \ \ \ \ \ \ \ \ Universit\'{e} de Lille 1}\\
{\small Case Western Reserve University \ \ \ \ \ \ \ \ \ \ \ \ \ UFR de Math\'{e}matiques }\\
{\small Cleveland, Ohio 44106, U. S. A. \ \ \ \ \ \ \ \ \ \ \ \ \ \ \ 59655 Villeneuve d'Ascq, France}\\
{\small \tt elisabeth.werner@case.edu}\\ \\


\vskip 5mm

\begin{thebibliography}{~~}
\small


\bibitem{GardnerBook}
{\sc R.J. Gardner}
{\em Geometric tomography.}
 Second edition. Encyclopedia of Mathematics and its Applications, 58. Cambridge University Press, Cambridge (2006).

\bibitem{Ga3}
{\sc R. J. Gardner}, {\em The dual Brunn-Minkowski theory for bounded Borel sets: Dual affine 
quermassintegrals and inequalities}, Adv. Math. {\bf 216} (2007), 358-386. 





\bibitem{Gruenbaum1963}
{\sc B. Gr\"unbaum},
{\em Measures of symmetry for convex sets},
Proc. Sympos. Pure Math. 7, (1963), 233--270.


\bibitem{Hab}
{\sc C. Haberl}, {\em Blaschke valuations}, Amer. J.  Math., 133, (2011), 717--751.











\bibitem{Lud2}
{\sc M. Ludwig}, {\em Ellipsoids and matrix valued valuations}, Duke Math. J. {\bf 119} (2003), 159-188.


\bibitem{Lud3}
{\sc M. Ludwig}, {\em Minkowski areas and valuations}, 
J. Differential Geometry, 86, (2010), 133--162.




\bibitem{LR2}
{\sc M. Ludwig and M. Reitzner,} {\em A classification of $SL(n)$
invariant valuations.}  Annals of Math. {\bf 172 } (2010), 1223-1271. 




\bibitem{Lu2}{\sc E. Lutwak}, {\em The Brunn-Minkowski-Firey theory II : Affine and geominimal surface areas}, Adv. Math. {\bf 118}  (1996),   244-294.



\bibitem{LYZ2002/1} {\sc E. Lutwak, D. Yang and G. Zhang}, {\em The Cramer--Rao inequality for star bodies}, Duke Math. J. {\bf 112} (2002), 59-81.

\bibitem{LYZ2004} {\sc E. Lutwak, D. Yang and G. Zhang}, {\em Volume inequalities for subspaces of $L_p$},   J. Differential Geom. {\bf 68} (2004), 159--184.




\bibitem{MeyerWerner1998} {\sc M. Meyer and E. Werner}, {\em The Santalo-regions of a convex body},
Transactions of the AMS {\bf 350}   (1998), 4569--4591.



\bibitem{MeyerSchuettWerner2011} 
{\sc M. Meyer, C. Sch\"utt and E. Werner},
{\em New affine measures of symmetry for convex bodies}, 
Adv. Math. {\bf 228}, (2011), 2920--2942.

\bibitem{MeyerSchuettWerner2012} 
{\sc M. Meyer, C. Sch\"utt and E. Werner}, 
{\em Affine invariant points}, 
{\em Preprint} 2012




\bibitem{NPRZ2010}
{\sc F. Nazarov,  F. Petrov, D. Ryabogin and A. Zvavitch},  {\em A remark on the Mahler conjecture: local minimality of the unit cube}, 
Duke Mathematical J. {\bf 154}, (2010), 419--430.

\bibitem{PaourisWerner2010} {\sc Relative entropy of cone measures and $L_p$ centroid bodies}
 Proc. London Math. Soc. (3) {\bf 104} (2012), no. 2, 253--286.


\bibitem{Rudin}
{\sc W. Rudin},
{\em Functional Analysis},
McGraw-Hill, 1991.



\bibitem{Santalo}
{\sc L.\ A.\ Santal\'o}, {\em Un invariante afin para los cuerpos convexos del
espacio de $n$ dimensiones\/}, Portugal.\ Math.\ {\bf 8} (1949), 155--161.



\bibitem{SchneiderBuch} {\sc R. Schneider}, 
{\em Convex Bodies: The
Brunn-Minkowski Theory}, Encyclopedia of Mathematics and its
Applications  44, Cambridge University Press, Cambridge (1993).



\bibitem{Schuster2010}
{\sc F. Schuster}, {\em Crofton measures and Minkowski valuations}, Duke Math. J. {\bf 154} (2010), 1--30.


\bibitem{SchuettWerner1990} {\sc C. Sch\"utt and E. Werner}, {\em The convex floating body},   {Math. Scand.} {\bf 66}, (1990), 275--290.



\bibitem{SW5}
{\sc C. Sch{\"u}tt and E. Werner}, {\em Surface bodies and
p-affine surface area.} Adv. Math. {\bf 187}  (2004), 98-145.



\bibitem{SA2}
{\sc A. Stancu}, {\em On the number of solutions to the
discrete two-dimensional $L_0$-Minkowski problem.} Adv. Math. {\bf
180} (2003), 290-323.


\bibitem{NicoleBuch} {\sc N. Tomczak Jaegerman}, 
{\em Banach-Mazur Distances and Finite-Dimensional Operator Ideals}, Encyclopedia of Mathematics and its
Applications  44, Cambridge University Press, Cambridge (1993).



 \bibitem{Werner1994} {\sc E. Werner}, {\em  Illumination bodies and affine surface area}, Studia
Math. {\bf 110}  (1994), 257--269.


 \bibitem{Werner2012} {\sc E. Werner}, {\em R\'enyi Divergence and $L_p$-affine surface area  for convex bodies}, 
Advances in Mathematics {\bf 230}, (2012), 1040--1059


\bibitem{WernerYe2008} {\sc E. Werner and D. Ye}, {\em New $L_{p}$ affine isoperimetric inequalities}, Adv. Math. {\bf 218} (2008), no. 3, 762-780.


\bibitem{WernerYe2010} {\sc E. Werner and  D. Ye}, {\em Inequalities for mixed $p$-affine surface area}, {Math. Ann.} {\bf  347}  (2010), 703-737.






\end{thebibliography}
\end{document}